\documentclass[12pt]{amsart}
\usepackage{amsmath,amscd,amssymb,amsfonts}
\newcommand{\h}{\hbox}
\newcommand{\q}{\quad}
\newcommand{\nin}{\noindent}
\newcommand{\bs}{\par\bigskip}
\newcommand{\ms}{\par\medskip}
\newcommand{\sk}{\par\smallskip}
\newcommand{\mprod}{\h{$\prod$}}
\newcommand{\msum}{\h{$\sum$}}
\newcommand{\mopl}{\h{$\bigoplus$}}
\newcommand{\ssb}{\raise.15ex\h{${\scriptscriptstyle\bullet}$}}
\newcommand{\C}{{\mathbb C}}
\newcommand{\DD}{{\mathbb D}}
\newcommand{\N}{{\mathbb N}}
\newcommand{\Q}{{\mathbb Q}}
\newcommand{\R}{{\mathbb R}}
\newcommand{\Z}{{\mathbb Z}}
\newcommand{\sbb}{{\mathbf s}}
\newcommand{\CC}{{\mathcal C}}
\newcommand{\D}{{\mathcal D}}
\newcommand{\E}{{\mathcal E}}
\newcommand{\OO}{{\mathcal O}}
\newcommand{\pr}{{\rm pr}}
\newcommand{\prh}{\widehat{\rm pr}}
\newcommand{\sih}{\widehat{\si}}
\newcommand{\Sh}{\widehat{S}}
\newcommand{\Y}{\widehat{Y}}
\newcommand{\Ch}{\widehat{C}}
\newcommand{\Hf}{H''_{\!f}}
\newcommand{\HfS}{H''_{\!f,S}}
\newcommand{\HFS}{H''_{\!F,S}}
\newcommand{\Hh}{\widehat{H}''}
\newcommand{\Hhf}{\widehat{H}''_{\!f}}
\newcommand{\HhFS}{\widehat{H}''_{\!F,S}}
\newcommand{\HhFSh}{\widehat{H}''_{\!F,\Sh}}
\newcommand{\HhfSh}{\widehat{H}''_{\!f,\Sh}}
\newcommand{\Gh}{\widehat{G}}
\newcommand{\Ghf}{\widehat{G}_{\!f}}
\newcommand{\GhfS}{\widehat{G}_{\!f,S}}
\newcommand{\GhFS}{\widehat{G}_{\!F,S}}
\newcommand{\GhfSh}{\widehat{G}_{\!f,\Sh}}
\newcommand{\GhFSh}{\widehat{G}_{\!F,\Sh}}
\newcommand{\So}{{}\,\overline{\!S}{}}
\newcommand{\ze}{\zeta}
\newcommand{\zo}{\overline{\zeta}}
\newcommand{\vv}{\overline{v}}
\newcommand{\et}{\widetilde{e}}
\newcommand{\st}{\widetilde{s}}
\newcommand{\Rt}{\widetilde{R}}
\newcommand{\tht}{\widetilde{\theta}}
\newcommand{\wt}{\widetilde{w}'}
\newcommand{\al}{\alpha}
\newcommand{\be}{\beta}
\newcommand{\ep}{\varepsilon}
\newcommand{\de}{{\delta}}
\newcommand{\De}{{\Delta}}
\newcommand{\ga}{\gamma}
\newcommand{\Ga}{\Gamma}
\newcommand{\la}{\lambda}
\newcommand{\lao}{\overline{\lambda}}
\newcommand{\dd}{{\partial}}
\newcommand{\dt}{\dd_t}
\newcommand{\dti}{\dd_t^{-1}}

\newcommand{\mo}{{\mathfrak m}_0}
\newcommand{\om}{\omega}
\newcommand{\Om}{\Omega}
\newcommand{\si}{\sigma}
\newcommand{\Gr}{{\rm Gr}}
\newcommand{\End}{{\rm End}}
\newcommand{\Hom}{{\rm Hom}}
\newcommand{\Imm}{{\rm Im}}
\newcommand{\bl}{\bigl}
\newcommand{\br}{\bigr}
\newcommand{\into}{\hookrightarrow}
\newcommand{\onto}{\mathop{\rlap{$\to$}\hskip2pt\h{$\to$}}}
\newcommand{\simto}{\buildrel\sim\over\longrightarrow}
\newcommand{\ges}{\geqslant}
\newcommand{\les}{\leqslant}
\newcommand{\plim}{\rlap{\raise-5.5pt\h{$\,\leftarrow$}}{\rm lim}}

\newcommand{\scd}{\h{$\cdot$}}
\newcommand{\Fls}{\raise-4pt\h{$\scriptstyle F$}}
\newcommand{\pl}{\,{+}\,}
\begin{document}
\title[Brieskorn lattices]
{On the structure of Brieskorn lattices, II}
\author[M. Saito]{Morihiko Saito}
\address{RIMS Kyoto University, Kyoto 606-8502 Japan}
\dedicatory{To the memory of Egbert Brieskorn}
\begin{abstract}
We give a simple proof of the uniqueness of extensions of good sections for formal Brieskorn lattices, which can be used in a paper of C. Li, S. Li, and K. Saito for the proof of convergence in the non-quasihomogeneous polynomial case. Our proof uses an exponential operator argument as in their paper, although we do not use polyvector fields nor smooth differential forms.
We also present an apparently simpler algorithm for an inductive calculation of the coefficients of primitive forms in the Brieskorn-Pham polynomial case.
In a previous paper on the structure of Brieskorn lattices, there were some points which were not yet very clear, and we give some explanations about these, e.g. on the existence and the uniqueness of primitive forms associated with good sections, where we present some rather interesting examples.
In Appendix we prove the uniqueness up to a nonzero constant multiple of the higher residue pairings in some formal setting which is different from the one in the main theorem. This is questioned by D. Shklyarov.
\end{abstract}
\maketitle
\centerline{\bf Introduction}
\bs\nin
Let $f:(X,0)\to(\De,0)$ be a holomorphic function on a complex manifold, where $\De$ is an open disk with coordinate $t$. Assume $X_0:=f^{-1}(0)$ has an isolated singularity at $0$. We have the associated Gauss-Manin system $G_f$ and the Brieskorn lattice $\Hf\subset G_f$, where $G_f$ is a regular holonomic $\D_{\De,0}$-module on which the action of $\dt$ is bijective, and $\Hf$ is a finite submodule over $\C\{t\}$ and also over $\C\{\!\{\dti\}\!\}$ (the latter comes from the theory of microdifferential operators \cite{SKK}), see \cite{B}, \cite{P}, \cite{SS}, \cite{bl}, etc. There is a surjection
$${\rm pr}_0:\Hf\onto\Hf/\dti\Hf\cong\Om_f\,\bl(:=\Om^{n+1}_{X,0}/df\wedge\Om^n_{X,0}\cong\C\{x\}/(\dd f)\br),$$
where $(\dd f)\subset\C\{x\}$ is the Jacobian ideal generated by the partial derivatives $\dd_{x_i}f$ with $x=(x_0,\dots,x_n)$ a local coordinate system of $(X,0)$, and $n:=\dim X_0=\dim X-1$.
\sk
For a $\C$-linear section $\si_0$ of $\pr_0$, set $I_0:=\Imm\,\si_0$.
We say that $\si_0$ is {\it good} in this paper if
$$tI_0\subset I_0+\dti I_0,\q\h{i.e.}\q t\si_0=\si_0A_0+\dti\si_0A_1\q\bl(A_0,A_1\in\End_{\C}\bl(\Om_f\br)\br).
\leqno(0.1)$$
Let $V$ be the filtration of Kashiwara \cite{K} and Malgrange \cite{M1} on $G_f$ indexed decreasingly by $\Q$ so that the action of $\dt t-\al$ on $\Gr^{\al}_VG_f$ is nilpotent. It induces the filtration $V$ on $\Hf$ and $\Om_f$.
A good section is called {\it very good} in this paper if it is strictly compatible with $V$. (It is called {\it good} in \cite{bl}.)
In the weighted homogeneous polynomial case, every good section is very good (see Proposition~3.1 below) although it does not hold in general.
The eigenvalues of $A_1$, which are called the exponents associated with a good section, do not necessarily coincide with the usual exponents defined as in \cite{St} unless the section is very good (see Example~4.1 below).
Note that $A_1$ is not necessarily semisimple in general (see \cite{bl}).
This causes a certain problem when we have to take an eigenvector of $A_1$ which generates the Jacobian ring over $\C\{x\}$.
It is needed to construct a primitive form associated with a good section satisfying the orthogonality condition for the canonical pairing.
\ms
The existence of a very good section is proved in \cite{bl} by using Deligne's canonical splitting of the mixed Hodge structure \cite{D} (which is applied to the canonical mixed Hodge structure on the vanishing cohomology \cite{St}) together with the relation with the Brieskorn lattice as in \cite{SS}.
Note that very good sections correspond to opposite filtrations to the Hodge filtration on the vanishing cohomology which are stable by the action of $N:=\log T_u$ where $T_u$ is the unipotent part of the monodromy (see \cite[Theorem 3.6]{bl}).
In the weighted homogeneous polynomial case, $N$ vanishes and the existence of very good sections is trivial so that we do not need to use the above arguments at all.
The orthogonality condition for the higher residue pairings in \cite{SK1}, \cite{SK2} follows from the orthogonality of the corresponding splitting of the Hodge filtration with respect to the canonical self-pairing of the vanishing cohomology, since the pairings can be identified with this self-pairing, see \cite{bl}.
Using the extension argument as below, we can get a unique primitive form associated with a very good section satisfying the orthogonality condition, see Remark~3.7 below. However, the existence and the uniqueness of the associated primitive form do not hold in general unless a good section is very good, see Examples~4.3 and 4.4 below.
\sk
Let $F:Y\to\De$ be a deformation of $f$ with $Y=X\times S$, $S=\De^m$, and $F|_{X\times\{0\}}=f$.
Here we assume that the singular locus $C$ of $(F,pr):Y\to\De\times S$ is {\it proper} over $S$.
Then the calculation of the Gauss-Manin system and the Brieskorn lattice can be reduced to the case $C\cap(X\times\{0\})=\{0\}$ by shrinking $S$ and restricting to an open neighborhood of each connected component of $C$.
We have the Gauss-Manin system $G_{F,S}$ and the Brieskorn lattice $\HFS\subset G_{F,S}$, where $G_{F,S}$ is a regular holonomic $\D_{\De\times S,0}$-module on which the action of $\dt$ is bijective, and $\HFS$ is a finite submodule over $\C\{t,\sbb\}$ and also over $\C\{\sbb\}\{\!\{u\}\!\}$ (see (1.1.1) for the latter).
Here $u:=\dti$, and $\sbb=(s_1,\dots,s_m)$ is the coordinate system of $\De^m\subset\C^m$. Let $\mo\subset\C\{\sbb\}:=\C\{s_1,\dots,s_m\}$ be the maximal ideal generated by the $s_i$. There is a surjection
$$\pr_S:\HFS\onto\HFS/\dti\HFS\cong\Om_{F,S}\,\bl(:=\Om^{n+1}_{Y/S,0}\big/dF\wedge\Om^n_{Y/S,0}\br),$$
together with the canonical isomorphisms
$$G_{F,S}|_0=G_f,\q\HFS|_0=\Hf,\q\Om_{F,S}|_0=\Om_f,$$
where $G_{F,S}|_0:=G_{F,S}\big/\,\mo\,G_{F,S}$, etc.
For a $\C\{\sbb\}$-linear section $\si_S$ of $\pr_S$, set $I_S:=\Imm\,\si_S$. We say that $\si_S$ is {\it good} if
$$tI_S\subset I_S+\dti I_S,\q\dd_{s_i}I_S\subset I_S+\dt I_S.
\leqno(0.2)$$
It is shown by B.~Malgrange (see \cite{M2}, \cite{M3}) that any good section $\si_0$ of $\pr_0$ can be uniquely extended to a good $\C\{\sbb\}$-linear section $\si_S$ of $\pr_S$ by solving Birkhoff's Riemann-Hilbert problem in this case, see also \cite{SK2}, \cite{H}, \cite{Sab}, etc.
(Here the orthogonality condition for the higher residue pairings can be reduced easily to the case $S=pt$.)
\sk
We can also consider the formal Gauss-Manin system $\Ghf$ and the formal Brieskorn lattice $\Hhf$, which are free modules of rank $\mu$ over $\C((u))$ and $\C[[u]]$ respectively (where $u=\dti$). They can be obtained by taking the $u$-adic completion of $G_f$ and $\Hf$ as in \cite{exp}.
There is a natural projection
$$\prh_0:\Hhf\onto\Hhf/\dti\Hhf\cong\Om_f,$$
where the last isomorphism follows from the $u$-adic completion argument.
\sk
We also have the formal Gauss-Manin system $\GhFSh$ and the formal Brieskorn lattice $\HhFSh$, which are free modules of rank $\mu$ over $\C((u))[[\sbb]]$ and $\C[[u,\sbb]]:=\C[[u,s_1,\dots,s_m]]$ respectively.
There is a natural projection
$$\prh_{\Sh}:\HhFSh\onto\HhFSh/\dti\HhFSh\cong\Om_{F,\Sh},$$
where $\Om_{F,\Sh}$ is the $\mo$-adic completion of $\Om_{F,S}$ so that $\Om_{F,\Sh}:=\Om_{F,S}\otimes_{\C\{\sbb\}}\C[[\sbb]]$.
We can define the notion of good sections $\sih_0$, $\sih_{\Sh}$ in the same way as in the convergent case by using the analogues of conditions (0.1) and (0.2) where $I_0$ is defined by $\Imm\,\sih_0$, and $I_S$ is replaced by $I_{\Sh}:=\Imm\,\sih_{\Sh}$. We have the following.
\ms\nin
{\bf Theorem~1.} {\it Any good $\C$-linear section $\sih_0$ of $\prh_0$ satisfying $(0.1)$ can be extended uniquely to a good $\C[[\sbb]]$-linear section $\sih_{\Sh}$ of $\prh_{\Sh}$ satisfying $(0.2)$ with $I_S$ replaced by $I_{\Sh}:=\Imm\,\sih_{\Sh}$.}
\ms
In fact, this easily follows from an assertion which is irrelevant to the action of $t$, see Theorem~1.4 below.
Theorem~1 does not seem to be stated explicitly in \cite{LLS}, although it seems to be used there in an essential way for the proof of the coincidence with the Malgrange's construction \cite{M2}, \cite{M3}, which gives the convergence of their extensions of good sections.
Here it seems rather difficult to prove directly the convergent version of Theorem~1 by using the exponential operator argument without using Malgrange's result in the convergent case.
The advantage of this method seems to be that one can calculate step by step the coefficients of the Taylor expansion of primitive forms explicitly (see (2.3) below for a special case). However, it is not very clear how much
it is useful for the original purpose of the primitive form, i.e. the associated period mapping, since the radius of convergence, for instance, does not seem to be calculated easily. It might be rather difficult to expect it theoretically since the partial Fourier transformation is used in an essential way.
\sk
It seems that Theorem~1 is proved in \cite{LLS} provided that ``uniquely" is replaced by ``canonically" in the statement.
In a more recent version of it, they seem to show the uniqueness statement in terms of primitive forms together with a rather complicated proof in the weighted homogeneous case.
Actually Theorem~1 can be proved more easily as is shown in the proof of Theorem~1.4 below by using an exponential operator argument given in \cite{LLS}. However, the latter argument is a rather amazing one for many complex geometers and their paper is not necessarily easy to read for non-specialists of mathematical physics. So we present in this paper a possibly simpler proof without using polyvector fields nor $C^{\infty}$ differential forms and by using a hopefully more precise argument than \cite{LLS}.
\sk
As a corollary of the exponential operator argument, we also present an algorithm for an inductive calculation of the coefficients of primitive forms for Brieskorn-Pham polynomials, which seems simpler than the one in \cite{LLS} in case of these polynomials.
By using it, we can calculate the coefficients of the first few terms of the Taylor expansion of the primitive forms without computers in this case, see (2.3) below. (The argument in this paper cannot be applied to the situation of \cite{DoSa} where the Brieskorn lattices are stable by $\dti$, but the $V$-filtration is stable by $\dt$, instead of $\dti$, in their case.)
\sk
In Appendix we prove the uniqueness up to a nonzero constant multiple of the higher residue pairings in some formal setting which is different from the one in Theorem~1 because of the difference between $\C((u))[[\sbb]]$ and $\C[[\sbb]]((u))$. It is written to answer a question of Dmytro Shklyarov. This uniqueness does not hold for the formal Gauss-Manin systems as in Theorem~1 because of the isomorphism in Proposition~1.3 below which is obtained by using the exponential operator argument. This shows a clear difference between the two kinds of formal Gauss-Manin systems.
\ms
We thank C.~Hertling for useful comments about this paper, D.~Shklyarov for a good question which became a source of Appendix, and C.~Li for a good question that led us to a correct formulation of an algorithm for the inductive calculation of the coefficients of primitive forms.
This work is partially supported by Kakenhi 24540039.
\ms
In Section 1 we review formal Gauss-Manin systems and Brieskorn lattices, and explain an exponential operator argument as in \cite{LLS}.
In Section 2 we present an algorithm for an inductive computation of the coefficients of the Taylor expansion of primitive forms in the Brieskorn-Pham polynomial case, which is apparently simpler in this case than the one in \cite{LLS}.
In Section 3 we give some remarks related to good sections and very good sections in the sense of this paper.
In Section 4 we present some interesting examples.
In Appendix we show the uniqueness up to a nonzero constant multiple of the higher residue pairings in some formal setting.
\bs\bs
\vbox{\centerline{\bf 1. Formal Gauss-Manin systems and Brieskorn lattices}
\bs\nin
In this section we review formal Gauss-Manin systems and Brieskorn lattices, and explain an exponential operator argument as in \cite{LLS} without using polyvector fields nor $C^{\infty}$ differential forms, but using more precise arguments.}
\ms\nin
{\bf Notation~1.1.}
Let $f:X\to\De$, and $F:Y\to\De$ be as in the introduction, where $Y=X\times S$ with $S=\De^m$.
We have the microlocal Gauss-Manin system defined by
$$G_{F,S}:=H^{n+1}C_{\!F,Y}^{\ssb}\q\h{with}\q C_{\!F,Y}^{\ssb}:=\bl(\,\Om^{\ssb}_{Y/S,0}\{\!\{u\}\!\}[u^{-1}],\,ud-dF\wedge\,\br),$$
where $u=\dti$, and $n=\dim X-1$.
Here $\Om^{\ssb}_{Y/S,0}\{\!\{u\}\!\}$ can be defined by using
$$\C\{y\}\{\!\{u\}\!\}:=\bl\{\msum_{\nu,k}\,a_{\nu,k}\,y^{\nu}u^k\in\C[[y,u]]\,\big|\,\msum_{\nu,k}\,|a_{\nu,k}|\,r^{|\nu|+k}/k!<\infty\,\,\,(\exists\,r>0)\br\},
\leqno(1.1.1)$$
where $y=(y_0,\dots,y_{n+m})$ is a local coordinate system of $Y$ with $y^{\nu}:=\prod_iy_i^{\nu_i}$ and $|\nu|:=\sum_i\nu_i$ for $\nu=(\nu_0,\dots,\nu_{n+m})\in\N^{n+m+1}$.
\sk
The Brieskorn lattice is defined by
$$\HFS:=H^{n+1}C_{\!F,Y}^{(0),\ssb}\q\h{with}\q C_{\!F,Y}^{(0),\ssb}:=
\bl(\,\Om^{\ssb}_{Y/S,0}\{\!\{u\}\!\},\,ud-dF\wedge\,\br).$$
These are obtained by the microlocalization of the usual Gauss-Manin systems and Brieskorn lattices, see \cite{P}, \cite{bl}, etc.
(Note that $G_{F,S}$ and $\HFS$ are finite free modules of rank $\mu$ over $\C\{\sbb\}\{\!\{u\}\!\}[u^{-1}]$ and $\C\{\sbb\}\{\!\{u\}\!\}$ respectively although it is not used in this paper.)
\sk
The action of $\dd_{x_j}$, $\dd_{s_i}$ can be defined by using the canonical generator $\de(t-F)$ which is not explicitly written in $C_{\!F,Y}^{\ssb}$ to simplify the notation (see also \cite{bl}).
More precisely $\de(t-F)$ is a generator of an $\E$-module $\CC_F$ which is the microlocalization of a $\D$-module ${\mathcal B}_F$, and the latter is the direct image of $\OO_Y$ by the graph embedding of $F$ as a $\D$-module. Here $\E$ is the ring of microdifferential operators (see \cite{SKK}). This generator satisfies the relations
$$\aligned t\,\de(t-F)&=F\,\de(t-F),\\
\dd_{x_j}\de(t-F)&=-(\dd F/\dd{x_j})\,\dt\,\de(t-F),\\
\dd_{s_i}\de(t-F)&=-(\dd F/\dd{s_i})\,\dt\,\de(t-F).\endaligned
\leqno(1.1.2)$$
Note that the second relation is compatible with the differential $ud-dF\wedge$ of the complex $C_{\!F,Y}^{\ssb}$ (up to the multiplication by $u$), and the latter can be identified with the relative de Rham complex ${\rm DR}_{Y/S}(\CC_F)$ up to a shift of complex. These are compatible with the theory of Gauss-Manin connections on Brieskorn lattices as in \cite{Gre}.
\sk
We have the formal Gauss-Manin system defined by
$$\GhFS:=H^{n+1}\Ch_{\!F,Y}^{\ssb}\q\h{with}\q\Ch_{\!F,Y}^{\ssb}:=\bl(\,\Om^{\ssb}_{Y/S,0}((u)),\,ud-dF\wedge\,\br),$$
see also \cite{sasa}, etc.\ for the case $S=pt$.
It has the formal Brieskorn lattice defined by
$$\HhFS:=H^{n+1}\Ch_{\!F,Y}^{(0),\ssb}\q\h{with}\q\Ch_{\!F,Y}^{(0),\ssb}:=
\bl(\,\Om^{\ssb}_{Y/S,0}[[u]],\,ud-dF\wedge\,\br).$$
We also have the bi-formal Gauss-Manin system defined by
$$\GhFSh:=H^{n+1}\Ch_{\!F,\Y}^{\ssb}\q\h{with}\q\Ch_{\!F,\Y}^{\ssb}:=\bl(\,\Om^{\ssb}_{X,0}((u))[[\sbb]],\,ud-dF\wedge\,\br),$$
with $[[\sbb]]:=[[s_1,\dots,s_{\mu}]]$, and similarly for $\HhFSh$ and $\Ch_{\!F,\Y}^{(0),\ssb}$ with $((u))$ replaced by $[[u]]$.
\sk
We can define similarly
$$G_{f,S},\q\GhfS,\q\GhfSh,\q\HfS,\q\Hh_{\!f,S},\q\HhfSh,$$
by replacing $F$ with $f$ in the above definitions, where $f$ is viewed as a trivial deformation.
\sk
We also have $\Gh_f$, $\Hhf$ by replacing $\Om_{Y/S,0}^{\ssb}$ with $\Om_{X,0}^{\ssb}$ in the definition of $\GhfS$, $\Hh_{\!f,S}$. There are canonical isomorphisms
$$\GhFSh\big|_0=\Gh_f,\q\HhFSh\big|_0=\Hhf,
\leqno(1.1.3)$$
and similar isomorphisms with $F$ replaced by $f$.
Here we set for any $\C[[\sbb]]$-module $N$
$$N|_0:=N/\mo N=N\otimes_{\C[[\sbb]]}\C,
\leqno(1.1.4)$$
where $\mo$ is the maximal ideal of $\C[[\sbb]]$. We also have a canonical injection
$$\iota:\Gh_f\into\GhfSh.
\leqno(1.1.5)$$
\sk
There are natural isomorphisms
$$\Om_{F,\Sh}=\HhFSh/\dti\HhFSh,\q\Om_{f,\Sh}=\HhfSh/\dti\HhfSh,\q\Om_f=\Hhf/\dti\Hhf,$$
where $\Om_{F,\Sh}$, $\Om_f$ are as in the introduction, and $\Om_{f,\Sh}=\Om_f[[\sbb]]$.
We have the canonical isomorphisms
$$\Om_{F,\Sh}\big|_0=\Om_f,\q\Om_{f,\Sh}\big|_0=\Om_f.
\leqno(1.1.6)$$
\ms\nin
{\bf Proposition~1.2.} {\it With the above notation, $\GhFSh$ and $\HhFSh$ are finite free modules of rank $\mu$ over $\C((u))[[\sbb]]$ and $\C[[u,\sbb]]=\C[[u,s_1,\dots,s_m]]$ respectively, where $\mu$ is the Milnor number of $f$. We have a similar assertion with $F$ replaced by $f$.}
\ms\nin
{\it Proof.} It is enough to show the assertion for $F$ since the assertion for $f$ is the special case of a trivial deformation.
\sk
Let $U^{\ssb}$ be the $\mo$-adic filtration on $\Ch_{\!F,Y}^{\ssb}$, $\Ch_{\!F,\Y}^{\ssb}$, i.e.
$$U^k\,\Ch_{\!F,Y}^{\ssb}=\mo^k\,\Ch_{\!F,Y}^{\ssb},\,\,\,\h{etc.}$$
Then $\Ch_{\!F,\Y}^{\ssb}$ is the $\mo$-adic completion of $\Ch_{\!F,Y}^{\ssb}$ so that
$$\Ch_{\!F,\Y}^{\ssb}=\rlap{\raise-10pt\h{$\,\,\,\scriptstyle k$}}\plim\,\Ch_{\!F,\Y}^{\ssb}/\mo^k\,\Ch_{\!F,\Y}^{\ssb}=\rlap{\raise-10pt\h{$\,\,\,\scriptstyle k$}}\plim\,\Ch_{\!F,Y}^{\ssb}/\mo^k\,\Ch_{\!F,Y}^{\ssb}.
\leqno(1.2.1)$$
Moreover the filtration $U$ induces a strict filtration on the complexes, and the induced filtration $U$ on the cohomology groups coincides with the $\mo$-adic filtration on these $\C[[\sbb]]$-modules so that
$$\GhFSh=\rlap{\raise-10pt\h{$\,\,\,\scriptstyle k$}}\plim\,\GhFSh/\mo^k\,\GhFSh=\rlap{\raise-10pt\h{$\,\,\,\scriptstyle k$}}\plim\,\GhFS/\mo^k\,\GhFS,
\leqno(1.2.2)$$
(and a similar assertion holds for the corresponding Brieskorn lattices).
These are shown by an argument similar to \cite{mhp}, \cite{exp} using the acyclicity of the complexes $\Gr_U^k\,\Ch_{\!F,Y}^{\ssb}$ except for the highest degree together with the Mittag-Leffler condition \cite{Gro}. Here the acyclicity follows from the canonical isomorphisms
$$\Gr_U^0\Ch_{\!F,Y}^{\ssb}\otimes_{\C}\Gr_U^k\,\C[[\sbb]]\simto\Gr_U^k\,\Ch_{\!F,Y}^{\ssb}.
\leqno(1.2.3)$$
\sk
Taking the cohomology of the last isomorphism and using the strictness of the filtration $U$, we then get the isomorphisms
$$\Gr_U^0\GhFS\otimes_{\C}\Gr_U^k\,\C[[\sbb]]\simto\Gr_U^k\,\GhFS\,(=\Gr_U^k\,\GhFSh).
\leqno(1.2.4)$$
This implies that $\GhFSh$ is free of rank $\mu$ over $\C((u))[[\sbb]]$ since $\Gr_U^0\GhFS=\Gh_f$ is free of rank $\mu$ over $\C((u))$.
The argument is similar for $\HhFSh$.
This finishes the proof of Proposition~1.2.
\ms\nin
{\bf Proposition~1.3} (compare to \cite{LLS}). {\it We have the exponential operator
$$\Psi:=e^{(F-f)/u}:\GhfSh\to\GhFSh,
\leqno(1.3.1)$$
which is an isomorphism of finite free $\C((u))[[\sbb]]$-modules with inverse given by
$$\Phi:=e^{(f-F)/u}:\GhFSh\to\GhfSh.
\leqno(1.3.2)$$
Moreover, these are compatible with the actions of $t$ and $\dd_{s_i}$.}
\ms\nin
{\it Proof.} Since $F-f\in\mo\OO_{Y,0}$, we can verify that $\Psi$ and $\Phi$ induce $\C((u))[[\sbb]]$-linear morphisms between the complexes $\Ch_{\!F,\Y}^{\ssb}$ and $\Ch_{f,\Y}^{\ssb}$, and these are inverse of each other.
Moreover they are compatible with the actions of $t$ and $\dd_{s_i}$ which are defined by using (1.1.2).
(For $t$, set $v:=u^{-1}=\dt$, which gives the Fourier transform of $t$, i.e. $t$ is identified with $-\dd_v$.)
This finishes the proof of Proposition~1.3.
\ms\nin
{\bf Theorem~1.4.} {\it Let $\si_{\Sh}:\Om_{F,\Sh}\to\HhFSh$ be a $\C[[\sbb]]$-linear section of the canonical projection $p_{F,\Sh}:\HhFSh\to\Om_{F,\Sh}$ satisfying the condition
$$\dd_{s_i}I_{\Sh}\subset I_{\Sh}+u^{-1}I_{\Sh}\q\h{with}\q I_{\Sh}:=\Imm\,\si_{\Sh}.
\leqno(1.4.1)$$
Such a section of $p_{F,\Sh}$ is uniquely determined by $I_0:=I_{\Sh}\big|_0\subset\Gh_f$ so that}
$$I_{\Sh}=\HhFSh\cap\Psi\bl(\iota\bl(I_0[u^{-1}]\br)[[\sbb]]\br).
\leqno(1.4.2)$$
\ms\nin
{\it Proof.} By the isomorphism (1.3.1), the assertion is equivalent to
$$\Phi(I_{\Sh})=\Phi\bl(\HhFSh\br)\cap\iota\bl(I_0[u^{-1}]\br)[[\sbb]]\q\h{in}\q\GhfSh.
\leqno(1.4.3)$$
We will show the inclusion $\subset$ together with the assertion that the right-hand side of (1.4.3) is isomorphic to $\Phi\bl(\Om_{F,\Sh}\br)$ by the projection $\Phi(p_{F,\Sh})$ so that it also gives a section of $\Phi(p_{F,\Sh})$.
\sk
By Propositions~1.2 and 1.3, $\HhFSh$ and $\Phi(\HhFSh)$ are free $\C[[u,\sbb]]$-submodules of $\GhFSh$ and $\GhfSh$ respectively with rank $\mu$. We have moreover
$$\mo^k\,\Phi(\HhFSh)=\Phi(\HhFSh)\cap\mo^k\,\GhfSh,
\leqno(1.4.4)$$
i.e. the inclusion $\Phi(\HhFSh)\into\GhfSh$ is strictly compatible with the $\mo$-adic filtration. This follows from the injective morphism of short exact sequences
$$\begin{array}{cccccccccccc}
0&\to&\mo^k\,\Phi(\HhFSh)&\to&\Phi(\HhFSh)&\to&\Phi(\HhFSh)/\mo^k\,\Phi(\HhFSh)&\to&0\\
&&\cap&&\cap&&\cap\\
0&\to&\mo^k\,\GhfSh&\to&\GhfSh&\to&\GhfSh/\mo^k\,\GhfSh&\to&0\end{array}$$
Here the injectivity of the last vertical morphism is reduced to the case $k=1$ by using the graded quotients $\Gr_U^j$ of the $\mo$-adic filtration $U$ together with isomorphisms similar to (1.2.4) (which hold also for $\Phi(\HhFSh)$ since it is a finite free $\C[[u,s]]$-module).
\sk
Using again the graded quotients $\Gr_U^j$ together with (1.4.4) and isomorphism similar to (1.2.4), we then get
$$\GhfSh=\Phi(\HhFSh)\oplus\iota\bl(u^{-1}I_0[u^{-1}]\br)[[\sbb]],
\leqno(1.4.5)$$
since
$$\Gh_f=\Hhf\oplus u^{-1}I_0[u^{-1}]\q\h{and}\q\Phi(\HhFSh)/\mo\Phi(\HhFSh)=\Hhf.$$
By (1.4.5) we get the isomorphism between the right-hand side of (1.4.3) and $\Phi\bl(\Om_{F,\Sh}\br)$.
\sk
It now remains to show
$$\Phi(I_{\Sh})\subset\iota\bl(I_0[u^{-1}]\br)[[\sbb]].
\leqno(1.4.6)$$
But this follows immediately from condition (1.4.1). In fact, $\GhfSh$ is identified with $\Gh_f[[\sbb]]$ so that any element of $\GhfSh$ has a Taylor expansion in $s$, and moreover, the above identification and $\Phi$ are compatible with the iterated actions of the $\dd_{s_i}$ and also with the restriction to $\sbb=0$. This finishes the proof of Theorem~1.4.
\ms\nin
{\bf Remarks~1.5.} (i) Formal Gauss-Manin systems and formal Brieskorn lattices are treated also in \cite{LLS} where the use of polyvector fields does not seem to be quite essential for them.
\ms
(ii) The commutativity of the projective limit and the cohomology does not seem to be explained in \cite{LLS}. Here the Mittag-Leffler condition as in \cite{Gro} is usually needed.
This point is not completely trivial even if we have the acyclicity of the complex except for the top degree.
For instance, it is not quite clear whether any surjective morphism of projective systems induces a surjective morphism by passing to the projective limit, unless we know that the Mittag-Leffler condition is satisfied for the projective system defined by the kernel, see \cite{Gro}.
This might be applied to the surjection from the top term of the complex to the cohomology, where the strictness of the last differential is related.
\ms
(iii) The construction in \cite{LLS} is slightly different from the one in earlier papers \cite{SK1}, \cite{SK2}, where the deformation $F$ of $f$ was defined over a space of dimension $\mu-1$, instead of $\mu$, and the value of $F$ together with the natural projection is used in order to define a morphism to a space $S$ of dimension $\mu$.
Note also that one gets a formal Gauss-Manin system of $\mu+1$ variables in \cite{LLS}, where the relative critical locus $C$ is finite and {\it flat} over $S$, although the image of $C$ in $S$ is the discriminant locus in \cite{SK1}, \cite{SK2}, since $F$ is used for the morphism to $S$.
\ms
(iv) It seems to be quite difficult to prove the convergent version of Theorem~1. Even in case $f=x^a+y^b$ with $1/a+1/b<1/2$, for instance, the convergence of the image of a monomial in $x,y$ by $\Psi$ seems to be quite non-trivial. (Note that, even if we get a divergent power series by this, it does not contradict the result of Malgrange since the procedure of extending good sections is not so simple.)
Here the calculation seems easier for $\Phi$.
It may be possible to show the convergence in $\sbb$ for each fixed degree part for the variable $u$ provided that we take a standard representative of the versal deformation of $f$ (i.e. $F=f+\sum_ig_is_i$ with $g_i$ monomial generators of the Jacobian ring).
\ms
(v) It does not seem to be very clear what kind of argument is used for the proof of the coincidence of the new construction of the higher residue pairings in \cite{LLS} with the old one.
It could be shown, for instance, by using the uniqueness (up to a constant multiple) of the pairing in the versal unfolding case by generalizing an argument in \cite[2.7]{bl} about the duality of simple holonomic $\E$-modules to the $\widehat{\E}$-module case and using the compatibility with the base change by $\{0\}\into S$ for the one variable case.
Here it does not seem easy to conclude it only by using the coincidence after taking the graded quotients of the Hodge filtration, since an automorphism of a filtered Gauss-Manin system of one variable is not necessarily the identity even if it induces the identity by taking the graded quotients.
(Note that a non-degenerate pairing can be identified with an isomorphism with the dual up to a shift of filtration. If there are two non-degenerate pairings, then we can compose one isomorphism with the inverse of the other so that we get an automorphism.)
\ms
(vi) If polyvector fields are used in the theory of primitive forms as in \cite{LLS}, one may have to divide a representative of a primitive form by a holomorphic relative differential form of the highest degree $\Omega_{Z/S}$ in order to get a representative in the polyvector fields. In this case one might get a ``primitive function" rather than a primitive form (and this may be more natural for the product structure). In the simple singularity case, it is a constant function, and this seems always possible provided that one can take the relative differential form $\Omega_{Z/S}$ to be the primitive form in the usual sense.
\bs\bs
\vbox{\centerline{\bf 2. Some explicit calculations}
\bs\nin
In this section we present an algorithm for an inductive computation of the coefficients of the Taylor expansion of primitive forms in the Brieskorn-Pham polynomial case, which is apparently simpler in this case than the one in \cite{LLS}.}
\ms\nin
{\bf 2.1.~Primitive forms.} In the notation of the introduction, assume $F$ is a miniversal deformation of $f$ as in \cite{LLS} so that
$$\dim S=\mu\,(:=\dim\OO_{X,0}/(\dd f)).$$
Let $\si_0:\Om_f\into\Hf$ be a good section of $\pr_0:\Hf\to\Om_f$ in (0.1) satisfying
$$S_K(\om,\om')\subset\C\,u^{n+1}\q\h{for}\,\,\,\om,\om'\in\Imm\,\si_0.
\leqno(2.1.1)$$
Here $u:=\dti$, and we denote in this paper the higher residue pairings by
$$S_K:G_f\times G_f\to K:=\C\{\!\{u\}\!\}[u^{-1}].
\leqno(2.1.2)$$
Note that
$$S_K(\om,\om')\subset\C\{\!\{u\}\!\}\,u^{n+1}\q\h{for any}\,\,\,\om,\om'\in\Hf.
\leqno(2.1.3)$$
This implies a rather strong restriction on $\Hf$.
\sk
By Malgrange's theory on Birkhoff's Riemann-Hilbert problem (see \cite{M2}, \cite{M3}), any good section $\si_0$ of $\pr_0:\Hf\onto\Hf/\dti\Hf\cong\Om_f$ in (0.1) can be uniquely extended to a good $\C\{\sbb\}$-linear section
$$\si_S:\Om_{F,S}\into\HFS$$
of
$$\pr_S:\HFS\onto\HFS/\dti\HFS\cong\Om_{F,S},$$
as is explained in the introduction.
Moreover the good section $\si_0$ is uniquely lifted to a $\C$-linear morphism
$$\si^{\nabla}_S:\Om_f\into\HFS,$$
so that
$$\Imm\,\si^{\nabla}_S\subset\Imm\,\si_S,\q\dd_{s_j}(\Imm\,\si^{\nabla}_S)\subset\dt\,(\Imm\,\si_S).
\leqno(2.1.4)$$
In fact, the second condition of (0.2) in the introduction implies an integrable connection on $\Om_S$ (by considering the action of $\dd_{s_j}$ on $I_S$ mod $\dt I_S$), and $\si^{\nabla}_S$ is defined by using the flat sections of this connection so that only the component of the second term $\dt I_S$ in the second condition of (0.2) remains (see \cite{SK1}, \cite{SK2}). Thus the second condition of (2.1.4) holds. Here (2.1.1) is also extended to the case of $\si^{\nabla}_S$.
Note that, by the uniqueness of the extension in Theorem~1.4, these constructions are compatible with the formal completion and we have similarly $\si^{\nabla}_{\Sh}$, etc.
\sk
Assume there is $\zo_0\in\Om_f$ which is an eigenvector of $A_1$ in (0.1), and generates $\Om_f$ over $\C\{x\}$. Set
$$\ze_0:=\si_0(\zo_0)\in\Hf.$$
In the weighted homogeneous polynomial case, we have up to a nonzero constant multiple
$$\ze_0=[dx_0\wedge\cdots\wedge dx_n],
\leqno(2.1.5)$$
where $x_0,\dots,x_n$ are coordinates such that $\sum_iw_ix_i\dd_{x_i}f=f$ with $w_i\in\Q_{>0}$.
(This follows from Proposition~3.1 below.)
\sk
The primitive form $\ze_S$ associated with $\si_0$ and $\zo_0$ is then defined by
$$\ze_S:=\si^{\nabla}_S(\zo_0)\in\HFS.$$
Similarly the formal primitive form $\ze_{\Sh}$ associated with $\si_0$ and $\zo_0$ is defined by
$$\ze_{\Sh}:=\si^{\nabla}_{\Sh}(\zo_0)\in\HhFSh.$$
The latter coincides with the image of $\ze_S$ in $\HhFSh$ by Theorem~1.4 together with a remark after (2.1.4).
\ms\nin
{\bf 2.2.~Relation with the exponential operators $\Psi$ and $\Phi$.} In the notation of (2.1) and Proposition~1.3, the formal primitive form $\ze_{\Sh}$ is the unique element of $\HhFSh$ satisfying
$$\Phi(\ze_{\Sh})=\iota(\ze_0)\q{\rm mod}\,\,\,\iota\bl(u^{-1}I_0[u^{-1}]\br)[[\sbb]],
\leqno(2.2.1)$$
where $I_0:=\Imm\,\si_0$, $u:=\dti$, and $\iota$ is as in (1.1.5). In fact, the uniqueness of $\ze_{\Sh}$ follows from the direct sum decomposition (1.4.5), and (2.2.1) holds since
$$\Phi(\ze_{\Sh})\big|_0=\ze_{\Sh}\big|_0=\ze_0,\q\dd_{s_j}\Phi(\ze_{\Sh})=\Phi(\dd_{s_j}\ze_{\Sh})\in\iota\bl(u^{-1}I_0[u^{-1}]\br)[[\sbb]],$$
where the last assertion follows from the proof of Theorem~1.4 together with the second condition of (2.1.4).
\sk
This characterization of formal primitive forms is compatible with the construction in \cite{LLS}, since (2.2.1) is equivalent to
$$\ze_{\Sh}=\Psi\bl(\iota(\ze_0)\br)\q{\rm mod}\,\,\,\Psi\bl(\iota\bl(u^{-1}I_0[u^{-1}]\br)[[\sbb]]\br).
\leqno(2.2.2)$$
\ms\nin
{\bf 2.3.~Case of Brieskorn-Pham polynomials.} Assume
$$f:=\msum_{i=0}^n\,x_i^{m_i}\q(m_i\ges 2),$$
i.e. $f$ is a Brieskorn-Pham polynomial.
In this case we can calculate the first few terms of the coefficients of the Taylor expansion of $\ze_{\Sh}$ without using a computer program as follows.
\sk
Set
$$\Ga:=\N^{n+1}\cap\mprod_{i=0}^n\,[0,m_i-2],$$
so that
$$\#\Ga=\mprod_{i=0}^n\,(m_i-1)=\mu.$$
\sk\nin
We have the natural coordinates $s_{\nu}$ of $S=\C^{\mu}$ for $\nu=(\nu_1,\dots,\nu_n)\in\Ga$. We may assume
$$F=f+\msum_{\nu\in\Ga}\,g_{\nu}s_{\nu}\q\h{with}\q g_{\nu}=x^{\nu}:=\mprod_i\,x_i^{\nu_i}\,\,\,(\nu\in\Ga).
\leqno(2.3.1)$$
Moreover we have the canonical good section $\si_0$ such that
$$I_0\,(:=\Imm\,\si_0)=\msum_{\nu\in\Ga}\,\C\,[g_{\nu}\om_0]\subset\Hf\q\h{with}\q
\om_0:=dx_0\wedge\cdots\wedge dx_n.$$
\sk
In the Brieskorn-Pham polynomial case we have for any $\nu=(\nu_0,\dots,\nu_n)\in\N^{n+1}$
$$\dt\,[x^{\nu}\om_0]=\frac{\nu_i-m_i+1}{m_i}[x^{\nu}x_i^{-m_i}\om_0]\q\h{if}\,\,\,\nu_i\ges m_i-1.
\leqno(2.3.2)$$
This implies
$$[x^{\nu}\om_0]=0\,\,\,\h{in}\,\,\,\Hf\q\h{if}\,\,\,\nu_i+1\in m_i\N\,\,\,\h{for some}\,\,\,i.
\leqno(2.3.3)$$
(These become more complicated in the general weighted homogeneous polynomial case.)
\sk
Let $\ze_{S,k}$ be the image of $\ze_S$ in $\HFS/\mo^{k+1}\HFS$, where $\mo$ is the maximal ideal of $\OO_{S,0}$, and $k$ is a positive integer (which may be determined by the computational ability). Set
$$A_k:=\bl\{a=(a_{\nu})\in\N^{\Ga}\,\big|\,|a|\les k\br\}\q\h{with}\q|a|:=\msum_{\nu\in\Ga}\,a_{\nu}.$$
For $a=(a_{\nu})\in\N^{\Ga}$, define
$$p(a)=\bl(p(a)_0,\dots,p(a)_n\br)\in\N^{n+1}\q\h{by}\q p(a)_i:=\msum_{\nu\in\Ga}\,\nu_ia_{\nu},$$
so that
$$g^a:=\mprod_{\nu}\,g_{\nu}^{a_{\nu}}=\mprod_{i,\nu}\,x_i^{\nu_ia_{\nu}}=\mprod_i\,x_i^{p(a)_i}=:x^{p(a)}.$$
Define further
$$q(a)=\bl(q(a)_0,\dots,r(a)_n\br),\,\,\,r(a)=\bl(r(a)_0,\dots,r(a)_n\br)\,\,\,\,\h{in}\,\,\,\,\N^{n+1},$$
by the condition
$$p(a)_i=q(a)_i\,m_i+r(a)_i\q\h{with}\q0\les r(a)_i<m_i\q(\forall\,i\in[0,n]).$$
In particular, we have
$$q(a)_i=\biggl\lfloor\frac{p(a)_i}{m_i}\biggr\rfloor.
\leqno(2.3.4)$$
(Note that $\lfloor\al\rfloor:=\max\{k\in\Z\mid k\les\al\}$ for $\al\in\R$.) Set
$$e_a=\msum_{i=0}^n\,q(a)_i\,-\,|a|,$$
and
$$A'_k:=\bl\{a\in A_k\,\,\big|\,\,e_a\ges 0,\,\,r(a)\in\Ga\br\}.$$
Note that the last condition $r(a)\in\Ga$ is equivalent to that $r(a)_i\ne m_i-1\,\,(\forall i)$.
\sk
Using the characterization of $\ze_{S,k}$ in (2.2.1), we then get the following Taylor expansion in $s$ by increasing induction on $|\nu|:=\sum_i\nu_i\les k\,$:
$$\ze_{S,k}=\sum_{a\in A'_k}\,c_a\,\dd_t^{-e_a}[\,g_{r(a)}\,s^a\,\om_0\,]\Fls\in\HFS/\mo^{k+1}\HFS,
\leqno(2.3.5)$$
with $c_a\in\C$, $s^a:=\mprod_{\nu\in\Ga}\,s_{\nu}^{a_{\nu}}$, and $g_{r(a)}=x^{r(a)}$ by definition. Here $[\eta]\Fls$ for $\eta\in\Om_{Y/S}^{n+1}$ denotes its class in $\HFS$ (mod $\mo^{k+1}$).
For $\om\in\Omega_{X,0}^{n+1}$, its class in $\Hf$ is simply denoted by $[\om]$. We have
$$[s^{\nu}\eta]\Fls=s^{\nu}[\eta]\Fls,$$
since the differential of the Gauss-Manin complex is $\OO_S$-linear. Note, however, that
$$[s^{\nu}\om]\Fls\ne s^{\nu}[\om]\q\bl(\h{i.e.,}\q[\om]\Fls\ne [\om]\br)\q\h{for}\q\om\in\Omega_{X,0}^{n+1}.$$
In fact, they belong to different groups $\HFS$ and $\HfS$ or $\Hf$. (This is related with a question of C.~Li. It is a source of an error in a previous version where the formula was too much simplified.)
\sk
By the characterization (2.2.1) the summation in (2.3.5) is actually taken over
$$A''_k:=\bl\{a\in A_k\,\,\big|\,\,\dt^{|a|}[g^a\om_0]\notin\dt I_0[\dt]\br\}.$$
In the Brieskorn-Pham polynomial case we have
$$\dt^{|a|}[g^a\om_0]\notin\dt I_0[\dt]\iff\dt^{|a|}[g^a\om_0]\in\Hf\setminus\{0\},
\leqno(2.3.6)$$
by (2.3.2) and (2.3.3). Using the last two formulas again, we then get
$$A'_k=A''_k,$$
together with the Taylor expansion (2.3.5) inductively.
\sk
The coefficients $c_a$ for $a\in A'_k$ are inductively determined by comparing the coefficients of both sides of (2.2.1). Since
$$e^{(f-F)\dd_t}=e^{-\sum_{\nu\in\Ga}g_{\nu}s_{\nu}\dd_t}=\mprod_{\nu\in\Ga}\,e^{-g_{\nu}s_{\nu}\dd_t},
\leqno(2.3.7)$$
we get by using (2.3.2)
$$c_a=-\sum_{0\les b< a}\,\Biggl((-1)^{|a-b|}\,\frac{c_{b}}{\,(a-b)!\,}\,\prod_{i=0}^n\,\prod_{k_i=1}^{q(a,b)_i}\frac{\,r(b)_i+p(a-b)_i-k_im_i+1\,}{m_i}\Biggr),
\leqno(2.3.8)$$
with
$$q(a,b)_i:=\biggl\lfloor\frac{r(b)_i+p(a-b)_i}{m_i}\biggr\rfloor.$$
\sk\nin
Here $\lfloor*\rfloor$ is as in a remark after (2.3.4), $(a-b)!:=\mprod_{\nu\in\Ga}(a_{\nu}-b_{\nu})!$, and we have by definition
$$b\les a\iff b_{\nu}\les a_{\nu}\,\,(\forall\,\nu\in\Ga),\q\h{and}\q b<a\iff b\les a\,\,\,\h{and}\,\,\,b\ne a.$$
\ms\nin
{\bf 2.4.~Example.} Assume $f=x_1^7+x_2^3$ and $k=6$. Then the $s^a=\prod_{\nu}s_{\nu}^{a_{\nu}}$ for $a\in A'_k\setminus\{0\}$ are
$$s_{(5,1)}^3,\q s_{(4,1)}s_{(5,1)}^2,\q s_{(5,1)}^6,\q s_{(4,1)}s_{(5,1)}^5,\q s_{(4,1)}^2s_{(5,1)}^4,\q s_{(3,1)}s_{(5,1)}^5.
\leqno(2.4.1)$$
The corresponding $g_{r(a)}=x^{r(a)}$ in (2.3.5) are respectively
$$x_1,\q 1,\q x_1^2,\q x_1,\q 1,\q 1,
\leqno(2.4.2)$$
and we have $e_a=0$ for $a\in A'_k$ in this case. We denote the corresponding coefficients $c_a$ by
$$c_{(1)},\,\,\,\dots\,\,\,,c_{(6)}.
\leqno(2.4.3)$$
Using (2.3.8), we first get
$$\aligned c_{(1)}&=\h{$\frac{1}{3!}\scd\frac{9\scd 2}{7^2\scd 3}=\frac{1}{7^2}$}\,,\\
c_{(2)}&=\h{$\frac{1}{2!}\scd\frac{8}{7^2\scd 3}=\frac{2^2}{7^2\scd 3}$}\,,\endaligned$$
and then verify that $c_{(3)},\dots,c_{(6)}$ are respectively equal to
$$\aligned \h{$-\frac{1}{6!}\scd\frac{24\scd 17\scd 10\scd 3\scd 4}{7^4\scd 3^2}
\pl\frac{1}{3!}\scd\frac{1}{7^2}\scd\frac{10\scd 3}{7^2\scd 3}$}
&=-\h{$\frac{17\scd 2^2}{7^4\scd 3^2}\pl\frac{5}{7^4\scd 3}
=\frac{-68+15}{7^4\scd 3^2}=-\frac{53}{7^4\scd 3^2}$}\,,\\
\h{$-\frac{1}{5!}\scd\frac{23\scd 16\scd 9\scd 2\scd 4}{7^4\scd 3^2}
\pl\frac{1}{2!}\scd\frac{1}{7^2}\scd\frac{9\scd 2}{7^2\scd 3}
\pl\frac{1}{3!}\scd\frac{2^2}{7^2\scd 3}\scd\frac{9\scd 2}{7^2\scd 3}$}
&=-\h{$\frac{23\scd 2^4}{7^4\scd 5\scd 3}\pl\frac{3}{7^4}\pl\frac{2^2}{7^4\scd 3}
=\frac{-368+45+20}{7^4\scd 5\scd 3}=-\frac{101}{7^4\scd 5}$}\,,\\
\h{$-\frac{1}{4!\scd 2!}\scd\frac{22\scd 15\scd 8\scd 4}{7^4\scd 3^2}
\pl\frac{1}{2!}\scd\frac{1}{7^2}\scd\frac{8}{7^2\scd 3}
\pl\frac{1}{2!}\scd\frac{2^2}{7^2\scd 3}\scd\frac{8}{7^2\scd 3}$}
&=-\h{$\frac{11\scd 5\scd 2^2}{7^4\scd 3^2}\pl\frac{2^2}{7^4\scd 3}\pl\frac{2^4}{7^4\scd 3^2}
=\frac{(-55+3+4)2^2}{7^4\scd 3^2}=-\frac{2^6}{7^4\scd 3}$}\,,\\
\h{$-\frac{1}{5!}\scd\frac{22\scd 15\scd 8\scd 4}{7^4\scd 3^2}
\pl\frac{1}{2!}\scd\frac{1}{7^2}\scd\frac{8}{7^2\scd 3}$}
&=-\h{$\frac{11\scd 2^3}{7^4\scd 3^2}\pl\frac{2^2}{7^4\scd 3}
=\frac{(-22+3)2^2}{7^4\scd 3^2}=-\frac{19\scd 2^2}{7^4\scd 3^2}$}\,.
\endaligned$$
The conclusion agrees with a calculation in \cite{LLS} using a different algorithm together with a computer program.
\bs\bs
\vbox{\centerline{\bf 3. Good sections and very good sections}
\bs\nin
In this section we give some remarks related to good sections and very good sections in the sense of this paper.}
\bs\nin
{\bf Proposition~3.1.} {\it In the notation of the introduction, any good section of $\pr_0$ is very good, if $f$ is a weighted homogeneous polynomial.}
\bs\nin
{\it Proof.} By definition (see (1.1.2)), $A_0$ in (0.1) is identified with the action of $f$ on the Jacobian ring $\C\{x\}/(\dd f)$, and it vanishes in the weighted homogeneous case.
Hence the image of the section is stable by the action of $\dd_tt$ which is identified with $A_1$. So the assertion follows.
\ms
The following proposition implies a formula for the dimension of the parameter space of very good sections satisfying the orthogonality condition for the self-duality in the case $N=0$ (including the weighted homogeneous polynomial case), see Corollary~(3.3) below.
\ms\nin
{\bf Proposition~3.2.} {\it Let $H$ be a finite dimensional $\C$-vector space with a finite filtration $F$. Let $S$ be a self-pairing of $H$ such that $S(F^pH,F^qH)=0$ for $p+q=m+1$, and the induced pairing of $\Gr_F^pH$ and $\Gr_F^qH$ is non-degenerate for $p+q=m$, where $m\in\Z$ is a fixed number. Assume $S$ is $(-1)^{m}$-symmetric, i.e. $S(u,v)=(-1)^{m}S(v,u)$. Set $e_p:=\dim\Gr^p_FH$. Then splittings $H=\bigoplus_kG^k$ of the filtration $F$ {\rm (}i.e. $F^PH=\bigoplus_{k\ges p}G^k)$ satisfying the condition $S(G^p,G^q)=0\,\,\,(p+q\ne m)$ are parametrized by $\C^{\,d(H,F,S)}$ with}
$$d(H,F,S):=\begin{cases}\msum_{p<q<m-p}\,e_pe_q+\msum_{p<m/2}\,\binom{e_p}{2}&\h{if $m$ is even,}\\
\msum_{p<q<m-p}\,e_pe_q+\msum_{p<m/2}\,\binom{e_p\,+\,1}{2}&\h{if $m$ is odd.}\end{cases}
\leqno(3.2.1)$$
\ms\nin
{\it Proof.} Let $\So$ denote the induced pairing of $\Gr^p_FH\times\Gr^{m-p}_FH$. We have $e_p=e_{m-p}$ since $\So$ is non-degenerate. Take bases $(\vv_{p,i})_{i\in[1,e_p]}$ of $\Gr^p_FH$ ($p\in\Z$) satisfying
$$\So(\vv_{p,i},\vv_{m-p,j})=\ep_p\,\de_{i,j}\q\h{with}\q\ep_p=\pm1,$$
where $\de_{i,j}=1$ if $i=j$, and 0 otherwise. Since $S(u,v)$ is $(-1)^{m}$-symmetric, we have
$$\ep_p=(-1)^{m}\ep_{m-p}.
\leqno(3.2.2)$$
We can lift $\vv_{p,i}$ to $v_{p,i}\in F^pH\subset H$ so that
$$S(v_{p,i},v_{q,j})=\ep_p\,\de_{p,m-q}\,\de_{i,j}\q\h{(with}\,\,\,\ep_p\,\,\,\h{as above)}.
\leqno(3.2.3)$$
This will be shown in Lemma~3.4 below.
(In the case of polarized Hodge structures as in the case of Corollary~(3.3) below, this easily follows from the Hodge decomposition.)
\sk
Set
$$I:=\bl\{(p,i)\in\Z^2\mid i\in[1,e_p]\br\},$$
where $[1,e_p]=\emptyset$ if $e_p=0$. Set
$$J:=\bl\{((p,i),(q,j))\in I^2\mid p<q\br\}\subset I^2.$$
Then any splitting of the filtration $F$ is expressed by
$$(\theta_{(p,i),(q,j)})\in\C^J,$$
since it defines a lift $w_{p,i}\in F^pH$ of $\vv_{p,i}\in\Gr_F^pH$ for each $(p,i)$ by
$$w_{p,i}:=v_{p,i}+\msum_{(q,j)\in I,\,q>p}\,\theta_{(p,i),(q,j)}\,v_{q,j}\in F^pH,$$
which is the image of $\vv_{p,i}$ by the splitting of the canonical surjection
$$F^pH\to\Gr_F^pH.$$
Note that the ambiguity of the splitting is given by the vector space
$${\rm Hom}(\Gr_F^pH,F^{p+1}H),
\leqno(3.2.4)$$
and its dimension is $\msum_{q>p}\,e_pe_q$ for each $p$.
\sk
The orthogonality condition of the splitting for the pairing $S$ is given by the relations
$$S(w_{p,i},w_{q,j})=0\q\h{for}\,\,\,((p,i),(q,j))\in R,$$
with
$$R:=\begin{cases}\bl\{((p,i),(q,j))\in I^2\mid p+q<m,\,(p,i)\les(q,j)\br\}&\h{if $m$ is even,}\\
\bl\{((p,i),(q,j))\in I^2\mid p+q<m,\,(p,i)<(q,j)\br\}&\h{if $m$ is odd.}\end{cases}$$
Here we use the lexicographic order on $I$, i.e. $(p,i)<(q,j)\iff p<q$ or $p=q,\,i<j$.
\sk
By (3.2.3) we have
$$S(w_{p,i},w_{q,j})=\begin{cases}
\,0&\h{if}\,\,\,p+q>m,\\
\,\ep_p\,\delta_{i,j}&\h{if}\,\,\,p+q=m,\end{cases}$$
and $S(w_{p,i},w_{q,j})$ for $p+q<m$ is given by
$$\aligned S(w_{p,i},w_{q,j})&=\ep_{m-q}\,\theta_{(p,i),(m-q,j)}+\ep_p\,\theta_{(q,j),(m-p,i)}\\
&+\msum_{(r,k)\in I,\,p<r<m-q}\,\ep_r\,\theta_{(p,i),(r,k)}\,\theta_{(q,j),(m-r,k)}.\endaligned
\leqno(3.2.5)$$
Here note that we have by (3.2.2)
$$\ep_{m-q}+\ep_p\ne 0\,\,\,\h{in the case where $(p,i)=(q,j)$ and $m$ is even.}
\leqno(3.2.6)$$
\sk
Consider the map
$$\ga:R\into J\q((p,i),(q,j))\mapsto((p,i),(m-q,j)).$$
We say that $\theta_{\ga((p,i),(q,j))}=\theta_{(p,i),(m-q,j)}$ is the {\it depending parameter} of the relation
$$S(w_{p,i},w_{q,j})=0\q\h{for}\q((p,i),(q,j))\in R.$$
By (3.2.5), $\theta_{(p,i),(m-q,j)}$ appears in $S(w_{p,i},w_{q,j})$ as a linear term with a nonzero coefficient, where (3.2.6) is used in the case $(p,i)=(q,j)$ and $m$ is even. Moreover $\theta_{(p',i'),(m-q',j')}$ appearing in the relation $S(w_{p,i},w_{q,j})=0$ must satisfy the inequality
$$p'+q'\ges p+q.$$
(In fact, $(p',i')$ must coincide with $(p,i)$ or $(q,j)$, and the inequality follows from (3.2.5).)
This implies that $\theta_{(p,i),(m-q,j)}$ does not appear in the relations
$$S(w_{p',i'},w_{q',j'})\q\h{with}\q p'+q'>p+q.$$
\sk
We can now prove by induction on $p+q$ and using (3.2.5) that the values of the depending parameters are given as polynomials of the remaining parameters
$$\theta_{(p,i),(q,j)}\q\h{with}\q((p,i),(q,j))\in J\setminus\ga(R),$$
which are called independent parameters. Thus splittings of the filtration $F$, which are orthogonal to each other with respect to the pairing $S$, are parametrized by
$$\C^{J\setminus\ga(R)}.$$
Moreover we have
$$d(H,F,S)=\#\bl(J\setminus\ga(R)\br).$$
So the assertion follows.
\ms\nin
{\bf Corollary~3.3} {\it Let $f:(X,0)\to(\De,0)$ be as in the introduction. Let $n=\dim X_0$. Assume the Milnor monodromy is semisimple. Let $n_{\al}$ be the multiplicity of the exponents of $f$ for $\al\in\Q\cap(0,n)$ as is defined in \cite{St}. Then very good sections of $\pr_0$ in the introduction are parametrized by $\C^{\,d_f}$ with $d_f=\sum_{|\la|=1,\,\Imm\,\la\ges 0}d_{f,\la}$ and
$$d_{f,\la}:=\begin{cases}\msum_{p<q<n+1-p}\,n_pn_q+\msum_{p<(n+1)/2}\,\binom{n_p}{2}&\h{if $\,\la=1$ and $n$ is odd,}\\
\msum_{p<q<n+1-p}\,n_pn_q+\msum_{p<(n+1)/2}\,\binom{n_p\,+\,1}{2}&\h{if $\,\la=1$ and $n$ is even.}\\
\msum_{p<q<n-p}\,n_{p+\al}n_{q+\al}+\msum_{p<n/2}\,\binom{n_{p+\al}}{2}&\h{if $\,\la=-1$ and $n$ is even,}\\
\msum_{p<q<n-p}\,n_{p+\al}n_{q+\al}+\msum_{p<n/2}\,\binom{n_{p+\al}\,+\,1}{2}&\h{if $\,\la=-1$ and $n$ is odd.}\\
\msum_{p<q}\,n_{p+\al}n_{q+\al}&\h{if $\,|\la|=1$ and $\,\Imm\,\la>0$,}\end{cases}$$
where $p,q\in\Z$, and $\la=e^{2\pi i\al}$ with $\al\in[0,\frac{1}{2}]$.}
\ms\nin
{\it Proof.} By \cite{St} there is a canonical mixed Hodge structure on the vanishing cohomology $H^n(F_{f,0},\C)$, where $F_{f,0}$ is the Milnor fiber of $f$ around $0\in X$, and the Hodge filtration $F$ is compatible with the direct sum decomposition by the eigenvalues of the monodromy $T$
$$H^n(F_{f,0},\C)=\mopl_{\la\in\C^*}\,H_{\la}.$$
Moreover there are canonical non-degenerate pairings of mixed Hodge structures
$$S:H_{\ne 1}\otimes H_{\ne 1}\to\C(-n),\q S:H_1\otimes H_1\to\C(-n-1),
\leqno(3.3.1)$$
where $H_{\ne 1}:=\mopl_{\la\ne 1}\,H_{\la}$, and these are compatible with the action of the monodromy $T$, i.e.
$$S(Tu,Tv)=S(u,v).
\leqno(3.3.2)$$
So the assumption on $S$ in Proposition~3.2 is satisfied for $H_{\ne 1}$ and $H_1$ with $m=n$ and $n+1$ respectively.
The multiplicities $n_{\al}$ of the Steenbrink exponents can be defined by
$$n_{\al}:=\dim\Gr_F^pH_{\la}\q\h{with}\q p=[\al],\,\,\la=e^{2\pi i\al},
\leqno(3.3.3)$$
where we use the symmetry of the exponents in \cite{St} i.e.
$$n_{\al}=n_{\be}\q\h{if}\q\al+\be=n+1.
\leqno(3.3.4)$$
\sk
For $\la=\pm 1$, the assertion of Corollary~(3.3) then follows from Proposition~3.2.
If $\la\ne\pm1$, we get the assertion by using the remark around (3.2.4) together with the duality isomorphism
$$\bl(H_{\lao},F[n]\br)={\bf D}(H_{\la},F):={\rm Hom}_{\C}\bl((H_{\la},F),\C\br),
\leqno(3.3.5)$$
which follows from the first non-degenerate pairing in (3.3.1).
(In fact, the latter implies that any splitting of $F$ on $H_{\la}$ determines uniquely its dual splitting of $F$ on $H_{\lao}$ by using the orthogonality condition with respect to $S$.)
This finishes the proof of Corollary~(3.3).
\ms\nin
{\bf Lemma~3.4} {\it With the notation in the proof of Proposition~$3.2$, the $\vv_{p,i}$ can be lifted to $v_{p,i}\in F^pH$ so that $(3.2.3)$ holds.}
\ms\nin
{\it Proof.} We show the assertion by induction on
$$\max\{p\mid\Gr^p_FH\ne 0\}-\min\{p\mid\Gr^p_FH\ne 0\}.$$
Set $a:=\min\{p\mid\Gr^p_FH\ne 0\}$, $b:=\max\{p\mid\Gr^p_FH\ne 0\}$, and
$$H'=F^{a+1}H/F^bH.$$
Let $S'$ be the induced pairing on $H'$.
By inductive hypothesis, $\vv_{p,i}$ for $p\in[a+1,b-1]$ can be lifted to $v'_{p,i}\in F^pH'\subset H'$ so that
$$S'(v'_{p,i},v'_{q,j})=\ep_p\,\de_{p,m-q}\,\de_{i,j}\q(p,q\in[a+1,b-1]).$$
We can lift $\vv_{a,i}$ to $v_{a,i}\in H$ by induction on $i$ so that
$$S(v_{a,i},v_{a,j})=0\q(i,j\in[1,e_a]).$$
Note that $\vv_{b,i}$ is identified with $v_{b,i}\in F^bH=\Gr^b_FH$, and we have
$$S(v_{a,i},v_{b,j})=\So(\vv_{a,i},\vv_{b,j})=\ep_a\,\de_{i,j}.$$
Then we can lift $v'_{p,i}$ to $v_{p,i}\in F^pH$ for $p\in[a+1,b-1]$ so that
$$S(v_{p,i},v_{a,j})=0\q(p\in[a+1,b-1]).$$
Here we have
$$S(v_{p,i},v_{q,j})=S'(v'_{p,i},v'_{q,j})=\ep_p\,\de_{p,m-q}\,\de_{i,j}\q(p,q\in[a+1,b-1]).$$
So (3.2.3) follows (since $S(v_{p,i},v_{b,j})=0$ for $p>a$).
This finishes the proof of Lemma~3.4.
\ms\nin
{\bf Remark~3.5.} In the weighted homogeneous polynomial case, it seems that the formula in Corollary~(3.3) is essentially equivalent to a formula for the parameter space of primitive forms in \cite{LLS}.
(Its verification is left to the reader.)
Condition (3.2.3) does not seem to be absolutely necessary for the argument in the proof of Proposition~3.2, since it seems to be enough to assume (3.2.3) for $p+q\ges m$ (which trivially holds) although (3.2.5) becomes more complicated without assuming condition (3.2.3) for $p+q<m$, see also \cite{LLS}.
Note, however, that the parameter space does not necessarily coincide with the origin in the case it is 0-dimensional, since it would imply (3.2.3) also for $p+q<m$.
\ms\nin
{\bf Remark~3.6.} We have in general
$$V^{>\al_{\mu}-1}\Hf=V^{>\al_{\mu}-1}G_f,
\leqno(3.6.1)$$
where $\al_{\mu}$ is the maximal exponent. In fact, setting $F^p\Hf:=\dt^{-p}\Hf$, we have
$$\Gr_F^p\Gr_V^{\al}\Hf=0\q\h{for}\q\al>\al_{\mu}+p,
\leqno(3.6.2)$$
(in particular, for $\al>\al_{\mu}-1$ and $p\les-1$).
\ms\nin
{\bf Remark~3.7.} It is known that the minimal exponent $\al_1$ in the usual sense (i.e. as is defined in (3.3.3)) has multiplicity 1, and moreover $V^{>\al_1}\Om_f\subset\Om_f$ is identified with the maximal ideal of the Jacobian ring $\C\{x\}/(\dd f)$, see \cite[4.11]{DiSa} (and also \cite{period}, Remark 3.11).
Here the theories of mixed Hodge modules \cite{mhp} and microlocal $b$-functions \cite{micro} are used. We need the commutativity of taking the graded quotients $\Gr_F^p$, $\Gr_V^{\al}$ and the cohomology functor $H^{n+1}$ in an essential way, since there is no canonical $\OO_X$-module structure if one takes the cohomology functor first.
(In case $\al_1<1$, the assertion may also follow from \cite{V}.)
\sk
The above assertion implies that there is a unique primitive form associated with any very good section (in the sense of this paper) satisfying the orthogonality condition for the higher residue pairings (which follows from the orthogonality condition as in \cite[Lemma~2.8]{bl}).
However, $A_1$ in (0.1) is not necessarily semisimple as is seen in Example~4.2 below, and there is not always a primitive form associated with any good section satisfying the orthogonality condition unless the section is very good, see Example~4.3 below.
We also have a problem about the uniqueness of the associated primitive form, see Example~4.4 below. If we assume that the eigenvalue of the Euler vector field is the minimal exponent, then this may make the existence of the associated primitive form more difficult in general.
\bs\bs
\vbox{\centerline{\bf 4. Examples.}
\bs\nin
In this section we present some interesting examples.}
\ms\nin
{\bf Example~4.1.} If $f$ is not a weighted homogeneous polynomial, it may be possible that there is a good section of $\pr_0$ which is not very good, see \cite{bl}. For instance, consider the case
$$f=x^a+y^b+x^{a-2}y^{b-2}\q(1/a+1/b<1/2),$$
where we have a good section such that the eigenvalues of $A_1$ in (0.1) are
$$\al'_1:=\al_1+1,\q\al'_{\mu}:=\al_{\mu}-1,\q\al'_k:=\al_k\,\,(k\in[2,\mu-1]).
\leqno(4.1.1)$$
Here $\al_1\les\cdots\les\al_{\mu}$ are the exponents of $f$ as is defined in \cite{St} (see also (3.3.3) above), which can be expressed in this case by
$$\msum_k\,t^{\,\al_k}=\msum_{0<i<a,\,0<j<b}\,t^{i/a+j/b},
\leqno(4.1.2)$$
with $\mu=(a-1)(b-1)$.
(Note that $\al'_i\les\al'_{i+1}$ does not hold for $i=1$ and $\mu-1$.)
\sk
To show (4.1.1), set
$$R:=\C\{\!\{\dti\}\!\},\q K:=\C\{\!\{\dti\}\!\}[\dt].
\leqno(4.1.3)$$
Put
$$\om^{(i,j)}=x^{i-1}y^{j-1}dx\wedge dy.$$
By using (1.1.2) restricted to $X\times\{0\}$, we get
$$\aligned t\,[\,\om^{(i,j)}]-\al^{(i,j)}\,\dti[\,\om^{(i,j)}]=c^{(i,j)}\,[\,\om^{(i+a-2,j+b-2)}]\q\h{in}\q\Hf,\\
\h{with}\q\q\al^{(i,j)}=\deg_{(a,b)}\om^{(i,j)}:=i/a+j/b,\q c^{(i,j)}\in\C^*.\endaligned
\leqno(4.1.4)$$
These imply that we have free generators $v_k$ ($k\in[1,\mu]$) of the Gauss-Manin system $G_f$ over $K$ satisfying
$$\dd_tt\,v_k=\al_kv_k\q(k\in[1,\mu]),
\leqno(4.1.5)$$
and we have the following free generators of the Brieskorn lattice $\Hf$ over $R:$
$$v_1+e\,\dt v_{\mu},\q v_k\,\,(k\in[2,\mu])\q\h{with}\q e\in\C^*.
\leqno(4.1.6)$$
More precisely the above calculation implies that
$$[\,\om^{(i,j)}]=v_k\,\,\,\,{\rm mod}\,\,\,\,V^{\al_k+2-2\al_1}G_f,
\leqno(4.1.7)$$
where $k$ is determined by $(i,j)\in[1,a-1]\times[1,b-1]$ with condition $i/a+j/b=\al_k$ satisfied.
Here $V$ is the filtration of Kashiwara and Malgrange on the Gauss-Manin system $G_f$ as in the introduction. This is closely related with the modified degree $\deg_{(a,b)}\om^{(i,j)}$ defined above, and we have
$$\deg_{(a,b)}\om^{(i,j)}\les\max\bl\{\al\in\Q\,\big|\,[\,\om^{(i,j)}]\in V^{\al}\Hf\br\},
\leqno(4.1.8)$$
where the equality holds if $(i,j)\in[1,a-1]\times[1,b-1]$. In fact, we have by \cite{exp}
$$\Gr_V^{\al_k}\om^{(i,j)}\ne 0\q\h{for}\,\,\,(i,j)\in[1,a-1]\times[1,b-1]\,\,\,\h{with}\,\,\,\al_k:=i/a+j/b.$$
(Here we can also use the $\mu$-constant deformation $f_s=x^a+y^b+s\,x^{a-2}y^{b-2}$ ($s\in\De^*$) together with the graded quotients of the decreasing filtration defined by $\deg_{(a,b)}\om\ges\al$ for $\om\in\Omega_X^2$.)
\sk
Take a good section whose image is spanned by
$$v'_1:=\dti v_1,\q v'_{\mu}:=\frac{1}{e}\,v_1+\dt\,v_{\mu},\q v'_k:=v_k\,\,(k\in[2,\mu-1]),
\leqno(4.1.9)$$
where $e\in\C^*$ is as above. Then the eigenvalues of the associated $A_1$ are as in (4.1.1).
\sk
Note that the image of $v'_{\mu}=\frac{1}{e}\,v_1+\dt\,v_{\mu}$ in the Jacobian ring modulo the maximal ideal does not vanish (i.e., it generates the Jacobian ring over it), and the other images vanish, where $\Om_X^2$ is trivialized by $dx\wedge dy$.
So $r$ in \cite{SK1}, \cite{SK2} seems to be $\al'_{\mu}=\al_{\mu}-1$ (instead of $\al_1$) which may be bigger than $\al_2$ in general.
It will be shown in Examples~4.3 and 4.4 below that this can cause serious problems related with the existence and the uniqueness of the associated primitive form.
\ms\nin
{\bf Example~4.2.} It is not very difficult to construct an {\it abstract} example of a Brieskorn lattice $\Hf$ with a good section such that $A_1$ in (0.1) is non-semi-simple.
(The following argument seems to be easier than the one in \cite{bl}, Remark after 3.10, where it seems rather difficult to determine the structure of the Brieskorn lattice for {\it geometric} examples.)
\sk
Let $(H',F)$ be the underlying filtered $\C$-vector space of a mixed $\R$-Hodge structure endowed with the self-duality pairing $S$, an automorphism $T_s$ of finite order, and a nilpotent endomorphism $N$ of type $(-1,-1)$, satisfying the usual conditions
$$S(T_su,T_sv)=S(u,v),\q S(Nu,v)+S(u,Nv)=0,\q T_sN=NT_s.$$
We have the eigenvalue decomposition $(H',F)=\mopl_{\la}\,(H'_{\la},F)$ by the action of $T_s$. Assume for simplicity
$$(H',F)=(H'_{\la},F)\oplus(H'_{\lao},F),$$
for some $\la\ne 1,-1$.
Then $(H'_{\lao},F)$ is the dual of $(H'_{\la},F)$ up to a shift of filtration by $S$. Assume further
$$\dim\Gr_F^pH'_{\la}=\begin{cases}1&\h{if}\,\,\,p=1\\2&\h{if}\,\,\,p=2,\\0&\h{otherwise,}\end{cases}\q\q
\dim\Gr_F^pH'_{\lao}=\begin{cases}2&\h{if}\,\,\,p=1,\\1&\h{if}\,\,\,p=2,\\0&\h{otherwise,}\end{cases}
\leqno(4.2.1)$$
together with the non-vanishing (i.e. the surjectivity and the injectivity) of the morphisms
$$N:\Gr_F^2H'_{\la}\onto\Gr_F^1H'_{\la},\q N:\Gr_F^2H'_{\lao}\into\Gr_F^1H'_{\lao}.$$
Then we have a splitting of the short exact sequence
$$0\to\Gr_F^2H'_{\la}\to H'_{\la}\to\Gr_F^1H'_{\la}\to 0,
\leqno(4.2.2)$$
such that the image of $\Gr_F^1H'_{\la}$ in $H'_{\la}$ by the splitting is contained in ${\rm Ker}\,N$, but does not coincide with $\Imm\,N$. For $H'_{\lao}$, we take the dual splitting by using $S$.
We will show that this splitting leads to an example of a good section of an {\it abstract} Brieskorn lattice $G_f^{\,\prime\,(0)}$ such that $A_1$ is non-semisimple.
\sk
By the above decompositions of $H'$, we have a decomposition of regular holonomic $\D_{S,0}$-modules
$$G'=G'_{\la}\oplus G'_{\lao}.
\leqno(4.2.3)$$
Here $G'$ is actually defined by the above isomorphism, and $G'_{\la}$, $G'_{\lao}$ are unique regular holonomic $\D_{S,0}$-modules of rank $3$ over $K$ together with isomorphisms
$$\Gr_V^{\be+k}G'_{\la}=H'_{\la},\q\Gr_V^{\be'+k}G'_{\lao}=H'_{\lao},
\leqno(4.2.4)$$
in a compatible way with the actions of $\dt t-\be-k$, $\dt t-\be'-k$, and $(2\pi i)^{-1}N$, where $\be,\be'\in\Q\cap(1,2)$ with $\la=e^{-2\pi i\be}$, $\lao=e^{-2\pi i\be'}$, and the action of $\dti$ is used for the above identification.
Then there are unique $R$-submodules $G_{\la}^{\,\prime\,(0)}$, $G_{\lao}^{\,\prime\,(0)}$ of $G'_{\la}$, $G'_{\lao}$ satisfying
$$\Gr_V^{\be+p}G_{\la}^{\,\prime\,(0)}=F^{2-p}H'_{\la},\q\Gr_V^{\be'+p}G_{\lao}^{\,\prime\,(0)}=F^{2-p}H'_{\lao}\q(\forall\,p\in\Z),
\leqno(4.2.5)$$
where $R,K$ are as in (4.1.3).
Moreover $G_{\la}^{\,\prime\,(0)}$ has free generators $e_1$, $e_2$, $e_3$ over $R$ satisfying
$$(\dt t-\be)\,e_1=\dt\,e_3,\q(\dt t-\be)\,e_2=0,\q(\dt t-\be-1)\,e_3=0.
\leqno(4.2.6)$$
(In fact, this follows from the vanishing of $\Gr_V^{\al}G_{\la}^{\,\prime\,(0)}$ for $\al\ne\be,\be+1$.)
\sk
The above choice of the splitting of (4.2.2) then gives free generators $\et_1$, $\et_2$, $\et_3$ of $G_{\la}^{\,\prime\,(0)}$ over $R$ defined by
$$\et_1:=e_1,\q\et_2:=e_2,\q\et_3:=e_3-c\,\dti e_2,
\leqno(4.2.7)$$
where $c\in\C^*$. Then we have
$$(\dt t-\be)\,\et_1=\dt\,\et_3+c\,\et_2,\q(\dt t-\be)\,\et_2=0,\q(\dt t-\be-1)\,\et_3=0.
\leqno(4.2.8)$$
So the action of $t$ on the generators $\et_1$, $\et_2$, $\et_3$ is expressed as in (0.1) by using the matrices
$$A_0=\begin{pmatrix}0&0&0\\ 0&0&0\\ 1&0&0\end{pmatrix}\q\q A_1=\begin{pmatrix}\be&0&0\\ c&\be&0\\0&0&\be+1\end{pmatrix}
\leqno(4.2.9)$$
and $A_1$ is non-semi-simple.
\ms\nin
{\bf Example~4.3.} It seems rather complicated to construct an example as in Example~4.2 above in a geometric way, and we need some more calculations as follows.
Here the Thom-Sebastiani type theorem as in \cite{SS} seems quite useful. For instance, set
$$f=g+h\q\h{with}\q g=x^{10}+y^3+x^2y^2,\,\,\,h=z^6+w^5+z^4w^3.$$
Let $G_f$, $\Hf$ denote the Gauss-Manin system and the Brieskorn lattice associated to $f$, and similarly with $f$ replaced by $g,h$.
Let $\al_{f,i}$ be the exponents of $f$, and similarly for $\al_{g,i}$, $\al_{h,i}$.
Then $H''_g$ has a basis $u_i$ over $R$ (with $R$ as in (4.1.3)) satisfying
$$(\dt t-\al_{g,1})\,u_i=\begin{cases}\dt\,u_{14}&\h{if}\,\,\,i=1,\\ 0&\h{if}\,\,\,i\ne 1,\end{cases}
\leqno(4.3.1)$$
where $\mu_g=14$, and we assume $\al_{g,i}\les\al_{g,i+1}$. In this case the $\al_{g,i}$ are given by
$$\msum_{i=1}^{14}\,t^{\al_{g,i}}=t^{1/2}+t+t^{3/2}+\msum_{k=1}^9\,t^{1/2+k/10}+\msum_{k=1}^2\,t^{1/2+k/3}.$$
In fact, this equality together with the non-triviality of the action of $N$ on $H_{-1}$ follows from a result in \cite{St} for functions with non-degenerate Newton boundary. Then (4.3.1) follows from \cite{SS} together with Remark~3.6, since
$$\al_{g,\mu_g}-\al_{g,1}=1.
\leqno(4.3.2)$$
\sk
As for $H''_h$, we have a basis $(v_1,\dots,v_{20})$ of $G_h$ over $K$ and free generators $v'_1,\dots,v'_{20}$ of $H''_h$ over $R$ satisfying (4.1.5) and (4.1.9) as in Example~4.1, where $\mu_h=20$, and $R,K$ are as in (4.1.3). We will denote $\al_j$, $\al'_j$ in (4.1.1) by $\al_{h,j}$, $\al'_{h,j}$ here.
\sk
We can actually take any $h$ in Example~4.1 satisfying the following condition:
$$\al_{g,i}+\al_{h,j}=\al_{g,\mu_g}+\al_{h,\mu_h}-2\q\h{for some}\,\,\,i,j\ges 2,
\leqno(4.3.3)$$
where $g$ may be replaced by $x^{a'}+y^{b'}+x^2y^2$ with $1/a'+1/b'<1/2$.
In the case of the above $g$ and $h$, condition (4.3.3) holds for $(i,j)=(2,2)$ as is shown later.
\sk
By the Thom-Sebastiani type theorem as in \cite{SS}, there are canonical isomorphisms
$$G_f=G_g\otimes_KG_h,\q\Hf=H''_g\otimes_RH''_h,
\leqno(4.3.4)$$
such that the action of $t$ on the left-hand side is identified with $t\otimes id+id\otimes t$ on the right-hand side.
Let $w_{i,j}$ and $w'_{i,j}$ be respectively the element of $G_f$ corresponding to $u_i\otimes v_j$ and $u_i\otimes v'_j$ in $G_g\otimes_KG_h$ under the isomorphism (4.3.4). Set
$$G'_f:=G'_{f,\la}\oplus G'_{f,\lao}\subset G_f,$$
with
$$\aligned G'_{f,\la}:&=K\,w_{1,20}\oplus K\,w_{2,2}\oplus K\,w_{14,20},\\
G'_{f,\lao}:&=K\,w_{1,1}\oplus K\,w_{13,19}\oplus K\,w_{14,1},\endaligned$$
where $\lambda=\exp(-2\pi i(2/15))$, and $\be=17/15$ in the notation of Example~4.2. In fact, we have
$$\begin{array}{llll}
\al_{g,1}=15/30,&\al_{g,2}=18/30,&\al_{g,13}=42/30,&\al_{g,14}=45/30,\\
\al_{h,1}=11/30,&\al_{h,2}=16/30,&\al_{h,19}=44/30,&\al_{h,20}=49/30,\end{array}$$
hence
$$\begin{array}{rrr}\al_{1,20}=32/15,&\al_{2,2}=17/15,&\al_{14,20}=47/15,\\
\al_{1,1}=13/15,&\al_{13,19}=43/15,&\al_{14,1}=28/15,\end{array}$$
and
$$\begin{array}{rrr}\al'_{1,20}=17/15,&\al'_{2,2}=17/15,&\al'_{14,20}=32/15,\\
\al'_{1,1}=28/15,&\al'_{13,19}=43/15,&\al'_{14,1}=43/15,\end{array}$$
where $\al_{i,j}:=\al_{g,i}+\al_{h,j}$, $\al'_{i,j}:=\al_{g,i}+\al'_{h,j}$.
Note that
$$(\dt t-\al_{i,j})^kw_{i,j}=0,$$
with $k=2$ if $i=1$, and $k=1$ otherwise.
\sk
If we consider the image of
$$R\,w'_{1,20}\oplus R\,w'_{2,2}\oplus R\,w'_{14,20},$$
by the natural projection $G'_f\onto G'_{f,\la}$, then it coincides with
$$R\,\dt w_{1,20}\oplus R\,w_{2,2}\oplus R\,\dt w_{14,20}.$$
So the situation is quite close to the one in Example~4.2.
\sk
Set
$$\wt_{i,j}:=\begin{cases}w'_{14,20}-c\,\dti w'_{2,2}&\h{if}\,\,\,(i,j)=(14,20)\\ w'_{13,19}+c'\dti w'_{1,1}&\h{if}\,\,\,(i,j)=(13,19)\\ w'_{i,j}&\h{otherwise.}\end{cases}$$
Here $c,c'\in\C^*$ are chosen appropriately so that $\wt_{14.20}$ and $\wt_{13,19}$ are orthogonal to each other.
Then $\wt_{i,j}$ and $\wt_{i',j'}$ are orthogonal to each other unless $(i,j)=(15-i',21-j')$.
Here we use the compatibility of the Thom-Sebastiani type isomorphism with the self-duality (i.e. with the higher residue pairings) up to a constant multiplication.
(This can be shown by using the fact that the discriminant of a deformation of the form $F:=f+\sum_ix_is_i$ is reduced.)
\sk
Let $G''_f$ be the orthogonal complement of $G'_f\subset G_f$ by the self-duality (i.e. the higher residue pairings). Then the decomposition $G_f=G'_f\oplus G''_f$ is compatible with the Brieskorn lattice, and induces the decomposition
$$\Hf=G_f^{\,\prime\,(0)}\oplus G_f^{\,\prime\prime\,(0)}.$$
In fact, we have the direct sum decompositions
$$\aligned G_g=G'_g\oplus G''_g\q\h{with}\q G'_g:=K\,u_1\oplus K\,u_{14},\q G''_g:=\mopl_{2\les i\les 13}\,K\,u_i,\\
G_h=G'_h\oplus G''_h\q\h{with}\q G'_h:=K\,v_1\oplus K\,v_{20},\q G''_h:=\mopl_{2\les i\les 19}\,K\,v_i,\endaligned$$
which are compatible with the Brieskorn lattices. They induce the decomposition compatible with the Brieskorn lattice
$$G_g\otimes_KG_h=(G'_g\otimes_KG'_h)\oplus(G''_g\otimes_KG''_h)\oplus(G'_g\otimes_KG''_h)\oplus(G''_g\otimes_KG'_h).$$
Then $G_f^{\,\prime\,(0)}$ is identified with the direct sum of
$$G'_g\otimes_KG'_h\q\h{and a direct factor of}\,\,\,G''_g\otimes_KG''_h,$$
via the isomorphism (4.3.4) in a compatible way with the Brieskorn lattice.
\sk
By a calculation similar to (4.2.8), the action of $t$ on the free generators
$$\wt_{1,20},\,\,\wt_{2,2},\,\,\wt_{14,20},\,\,\wt_{1,1},\,\,\wt_{13,19},\,\,\wt_{14,1}$$
of $G_f^{\,\prime\,(0)}$ over $R$ can be expressed as in (0.1) by using the matrices
$$A_0=\begin{pmatrix}0&0&0&0&0&0\\ 0&0&0&0&0&0\\ 1&0&0&0&0&0\\ \ga&0&0&0&0&0\\ 0&0&0&0&0&0\\ 0&0&\ga&1&0&0\end{pmatrix}\q\q A_1=\begin{pmatrix}\be&0&0&0&0&0\\ c&\be&0&0&0&0\\ 0&0&\be+1&0&0&0\\ 0&0&0&\be'&0&0\\ 0&0&0&0&\be'+1&0\\ 0&0&0&0&c'&\be'+1\end{pmatrix}
\leqno(4.3.5)$$
where $\be=17/15$, $\be'=28/15$, and $\ga\in\C^*$.
In this case it is rather difficult to get an associated primitive form. In fact, $\wt_{1,20}$ is the unique member of the generators whose class in the Jacobian ring $\OO_{X,0}/(\dd f)$ generates the ring over it, where $\Om_X^2$ is trivialized by $dx\wedge dy$. However, $\wt_{1,20}$ is annihilated only by $(A_1-\be)^2$, and the kernel of $A_1-\be$ in the Jacobian ring is generated over $\C$ by the class of $\wt_{2,2}=w_{2,2}$ which is contained in the maximal ideal.
(The details are left to the reader.)
\ms\nin
{\bf Example~4.4.} We first consider an abstract example.
Let $G$ be a regular holonomic $\D_{S,0}$-module which is a free $K$-module of rank 4 with generators $u_i\,(i\in[1,4])$ satisfying
$$\dt t\,u_i=\ga_iu_i,$$
with
$$0<\ga_1<\ga_k<\ga_4<1\q(k=2,3).
\leqno(4.4.1)$$
Assume $u_i$ and $u_j$ are orthogonal to each other by the self-duality pairing (i.e. the higher residue pairings) $S_K$ in (2.1.2) unless $i+j=5$. More precisely, assume
$$S_K(u_i,u_j)=\ep_i\,\de_{i,5-j}\,\dti,$$
with $\ep_i\in\C^*$ satisfying $\ep_1=\ep_2=-\ep_3=-\ep_4$. Note that the above condition implies $\ga_i+\ga_{5-i}=1$.
\sk
Let $c,c'\in\C^*$. Put
$$u'_i:=\begin{cases}u_1+cu_3+c'u_4&\h{if}\,\,\,i=1,\\ u_2+cu_4&\h{if}\,\,\,i=2,\\ \dti\,u_i&\h{if}\,\,\,i=3,4.
\end{cases}$$
Then
$$S_K(u'_i,u'_j)=\ep'_i\,\de_{i,5-j}\,\dt^{-2}\q(\ep'_i\in\C^*).$$
Set $c'':=c'/c$. Define
$$\aligned w_1:&=\dti u_1,\\ w_2:&=u'_1-c''u'_2=u_1-c''u_2+cu_3,\\ w_3:&=\dti u_3,\\ w_4:&=u'_1=u_1+cu_3+c'u_4.\endaligned$$
Then we have
$$\Hf:=\msum_{i=1}^4\,R\,u'_i=\msum_{i=1}^4\,R\,w_i,$$
and moreover
$$S_K(w_i,w_j)=\ep''_i\,\de_{i,5-j}\,\dt^{-2}\q(\ep''_i\in\C^*).$$
In this case the action of $t$ on the generators $w_1,\dots,w_4$ can be expressed as in (0.1) by using the matrices
$$A_0=\begin{pmatrix}0&*&0&*\\ 0&0&0&0\\ 0&*&0&*\\ 0&0&0&0\\ \end{pmatrix}\q\q A_1=\begin{pmatrix}\ga_1+1&0&0&0\\ 0&\ga_2&0&0\\ 0&0&\ga_3+1&0\\ 0&0&0&\ga_4\end{pmatrix}
\leqno(4.4.2)$$
\sk
This abstract example can be realized as a direct factor of the Brieskorn lattice associated with
$$f=x^a+y^b+x^{a-3}y^{b-2}+x^{a-2}y^{b-2},$$
if $a>b$ and $3/a+2/b<1$ (where the last condition corresponds to (4.4.1)). In fact, setting
$$g_1:=1,\q g_2:=x,\q g_3:=x^{a-3}y^{b-2},\q g_4:=x^{a-2}y^{b-2},$$
we have
$$\aligned{}[g_i\,dx\wedge dy]&=u_i\,\,\,\,{\rm mod}\,\,\,\,V^{\ga_i+2-\ga_1-\ga_2}G_f\q(i=1,2),\\
\dt\,[g_i\,dx\wedge dy]&=u_i\,\,\,\,{\rm mod}\,\,\,\,V^{\ga_i+1-\ga_1-\ga_2}G_f\q(i=3,4),\endaligned$$
where
$$\ga_1=1/a+1/b,\q\ga_2=2/a+1/b,\q\ga_3=1-2/a-1/b,\q\ga_4=1-1/a-1/b.$$
The argument is similar to the proof of (4.1.7). (The details are left to the reader.)
In this case, both $w_2$ and $w_4$ can be a primitive form associated with the good section whose image is spanned by the $w_i$.
\bs\bs
\vbox{\centerline{\bf Appendix: Uniqueness of higher residue pairings in some formal setting}
\bs\nin
This Appendix is written to answer a question of Dmytro Shklyarov.}
\ms\nin
Let $R=\C[[\sbb]]$ with $\sbb=(s_1,\dots,s_m)$, and $u:=\dti$.
Let $\Gh_R$ and $\Hh_R$ respectively denote the `formal' Gauss-Manin system and the `formal' Brieskorn lattice associated with a deformation $F=f+\sum_{i=1}^mg_is_i$ of $f\in\C\{x\}$ with an isolated singularity. Here `formal' means that $\Gh_R$ and $\Hh_R$ are finite free modules of rank $r$ over $R((u))$ and $R[[u]]$ respectively. They are endowed with the actions of $t$ and $\dd_{s_i}$ or $u\dd_{s_i}$ satisfying the usual relations. (Note that the uniqueness of the higher residue pairings does not hold over $\C((u))[[\sbb]]$ because of the isomorphism in Proposition~1.3. In fact, $\C((u))[[\sbb]]$ is much bigger than $R((u))$, and has much larger flexibility as is shown by the proposition.)
\sk
The dual of $\Gh_R$ can be defined by
$$\DD(\Gh_R):=\Hom_{R((u))}\bl(\Gh_R,R((u))\br),$$
where the actions of $R((u))$, $t$, and $\dd_{s_i}$ are given appropriately as usual, see e.g. \cite{bl}. Then the self-duality pairing (i.e. the higher residue pairings) can be identified with an isomorphism of $R((u))\langle\dd_{s_i},t\rangle$-modules
$$\Gh_R\simeq\DD(\Gh_R).$$
So the uniqueness up to a nonzero constant multiple of the higher residue pairings in this formal setting is equivalent to
$$\End_{R((u))\langle\dd_{s_i},t\rangle}(\Gh_R)=\C,
\leqno(A.1)$$
under the assumption that the discriminant is {\it reduced}, e.g. if $F$ is a miniversal deformation of $f$.
Here the discriminant $D$ is a divisor on $(\C\times\C^m,0)$ having the coordinates $t, s_1,\dots,s_m$, and $D$ is the image of the relative critical locus defined by the $\dd_{x_i}F$. We can also get $D$ by using the graded quotients of the filtration on the usual Gauss-Manin system defined by the usual Brieskorn lattice shifted by the action of $\dt^{-i}$, where the latter is a coherent sheaf on $(\C\times\C^m,0)$.
Passing to the completion by the maximal ideal of $\C\{\sbb\}$, we get the isomorphisms of $R[t]$-modules
$$\Hh_R/\dt^{-k}\Hh_R\cong R[t]/(h)^k,
\leqno(A.2)$$
where $h\in\C\{\sbb\}[t]$ is a defining function of the discriminant $D$, and $\Om_X^{n+1}$ is trivialized by $dx_0\wedge\cdots\wedge dx_n$.
\sk
There is a divisor $\Sigma$ on $(\C^m,0)$ such that $D\subset\C\times\C^m$ is etale over the complement of $\Sigma$ by the projection $\C\times\C^m\to\C^m$.
By Hironaka's resolution of singularities using blowing-ups with smooth centers, the assertion can be reduced to the case where $\Sigma$ is a divisor with normal crossings. In fact, the pull-back induces an injective morphism of local rings under smooth center blow-ups of $\C^m$, and we still have the injectivity after taking the formal completion for $s_i$ and $u$.
Then, changing the coordinates $s_i$ appropriately, we may assume that the discriminant $D$ is defined in $(\C\times\C^m,0)$ by the function
$$h:=t^r-s_1^{a_1}\cdots s_m^{a_m}.
\leqno(A.3)$$
Here we can forget the relation with $f,F$ from now on.
\sk
We take the ramified covering
$$\rho:(\C^m,0)\ni(\st_i)\mapsto(s_i):=(\st_i^{\,b_i})\in(\C^m,0),$$
where $b_i:=r/{\rm GCD}(r,a_i)$. Set $c_i:=a_i/{\rm GCD}(r,a_i)$. Then $rc_i=a_ib_i$, and the pull-back of the equation (A.3) under $\rho$ is given by
$$\widetilde{h}:=t^r-(\st_1^{\,c_1}\cdots\st_m^{\,c_m})^r.$$
\sk
We now pass to the localization $\Rt_{\st}:=\Rt[1/\st_1\cdots\st_m]$ of $\Rt:=\C[[\st_1,\dots,\st_m]]$. This is a finite etale Galois extension of $R_s:=R[1/s_1\cdots s_m]$ with Galois group $G=\prod_{i=1}^m\mu_{b_i}$, where $\mu_{b_i}$ is the group of roots of $1$ of order $b_i$ in $\C$. Let $\Gh_{\Rt_{\st}}$ be the pull-back of $\Gh_{R_s}:=R_s\otimes_R\Gh_R$ by $\rho$. This can be defined by $\Rt_{\st}\otimes_{R_s}\Gh_{R_s}$ since $\Rt_{\st}$ is finite over $R_s$. We have the canonical decomposition
$$\Gh_{\Rt_{\st}}=\mopl_{\lambda\in\mu_r}\,\Gh_{\Rt_{\st},\lambda},
\leqno(A.4)$$
where $\mu_r:=\{\lambda\in\C\mid\lambda^r=1\}$. In fact, let $F$ be the decreasing filtration on $\Gh_{\Rt_{\st}}$ defined by $u^j\Hh_{\Rt_{\st}}$ where $\Hh_{\Rt_{\st}}$ is the localization by $\st_1\cdots\st_n$ of the pull-back by $\rho$ of the formal Brieskorn lattice. Then we can get the decomposition by taking the inductive limit by $p$ of the projective limit by $q$ of the canonical decompositions
$$(F^p/F^q)\Gh_{\Rt_{\st}}=\mopl_{\lambda\in\mu_r}\,(F^p/F^q)\Gh_{\Rt_{\st},\lambda},
\leqno(A.5)$$
which can be defined by setting
$$(F^p/F^q)\Gh_{\Rt_{\st},\lambda}={\rm Ker}\bl((t-\lambda\,\st_1^{\,c_1}\cdots\st_m^{\,c_m})^{q-p}:(F^p/F^q)\Gh_{\Rt_{\st}}\to(F^p/F^q)\Gh_{\Rt_{\st}}\br),
\leqno(A.6)$$
since the discriminant is reduced.
In fact, there is a canonical direct sum decomposition
$$\aligned&\C\bl[t,\st_1,\dots,\st_m,\h{$\frac{1}{\st_1\cdots\st_m}$}\br]\big/\bl(t^r-(\st_1^{\,c_1}\cdots\st_m^{\,c_m})^r\br)^{q-p}\\&\q=\mopl_{\lambda\in\mu_r}\,\C\bl[t,\st_1,\dots,\st_m,\h{$\frac{1}{\st_1\cdots\st_m}$}\br]\big/\bl(t-\lambda\,\st_1^{\,c_1}\cdots\st_m^{\,c_m}\br)^{q-p}.\endaligned$$
Taking its tensor product with $\Rt=\C[[\st_1,\dots,\st_m]]$ over $\C[\st_1,\dots,\st_m]$, we then get
$$\aligned(F^p/F^q)\Gh_{\Rt_{\st}}&\cong\Rt\bl[t,\h{$\frac{1}{\st_1\cdots\st_m}$}\br]\big/\bl(t^r-(\st_1^{\,c_1}\cdots\st_m^{\,c_m})^r\br)^{q-p}\\&=\mopl_{\lambda\in\mu_r}\,\Rt\bl[t,\h{$\frac{1}{\st_1\cdots\st_m}$}\br]\big/\bl(t-\lambda\,\st_1^{\,c_1}\cdots\st_m^{\,c_m}\br)^{q-p},\endaligned
\leqno(A.7)$$
where the first isomorphism follows from (A.2). (In fact, $\rho$ is flat and the pull-back is an exact functor.)
This implies that the decomposition (A.5) can be obtained by (A.6).
(Note that $F$ cannot be exhaustive if we use the formal Gauss-Manin system as in Theorem~1.)
\sk
For $\theta\in\End_{R((u))\langle\dd_{s_i},t\rangle}(\Gh_R)$, its pull-back $\tht:=\rho^*\theta$ is an endomorphism of $\Gh_{\Rt_{\st}}$ preserving the decomposition (A.4). (In fact, $\tht$ preserves the filtration $F$ up to a shift by some integer $k$, i.e., $\tht(F^p\Gh_{\Rt_{\st}})\subset F^{p-k}\Gh_{\Rt_{\st}}$ for any $p$.) Moreover $\tht$ is compatible with the action of $G$ (since it is the pull-back of $\theta$ by $\rho$), and $G$ acts on the direct factors of the decomposition (A.4) transitively. Thus the assertion is reduced to
$$\End_{\Rt_{\st}((u))\langle\dd_{\st_i},t\rangle}(\Gh_{\Rt_{\st},\lambda})=\C.
\leqno(A.8)$$
We can verify (A.8) easily since $\Gh_{\Rt_{\st},\lambda}$ is a free $\Rt((u))\bl[\frac{1}{\st_1\cdots\st_m}\br]$-module of rank 1 by (A.7).
So (A.1) follows.


\begin{thebibliography}{DoSa}
\bibitem[Br]{B} Brieskorn, E., Die Monodromie der isolierten Singularit\"aten von Hyperfl\"achen, Manuscripta Math., 2 (1970), 103--161.
\bibitem[De]{D} Deligne, P., Th\'eorie de Hodge, II, Inst.\ Hautes Etudes Sci.\ Publ.\ Math.\ 40 (1971), 5--57.
\bibitem[DiSa]{DiSa} Dimca, A.\ and Saito, M., Koszul complexes and spectra of projective hypersurfaces with isolated singularities, arXiv:1212.1081.
\bibitem[DoSa]{DoSa} Douai, A.\ and Sabbah, C., Gauss-Manin systems, Brieskorn lattices and Frobenius structures, I, Ann.\ Inst.\ Fourier 53 (2003), 1055--1116.
\bibitem[Gre]{Gre} Greuel, G.-M., Der Gauss-Manin-Zusammenhang isolierter Singularit\"aten von vollst\"andigen Durchschnitten, Math.\ Ann.\ 214 (1975), 235--266.
\bibitem[Gro]{Gro} Grothendieck, A., El\'ements de g\'eom\'etrie alg\'ebrique, III-1, Publ.\ Math.\ IHES 11, 1961.
\bibitem[He]{H} Hertling, C., Frobenius Manifolds and Moduli Spaces for Singularities, Cambridge University Press, 2002.
\bibitem[Ka]{K} Kashiwara, M., Vanishing cycle sheaves and holonomic systems of differential equations, Lect. Notes in Math.\ 1016, Springer, Berlin, 1983, pp.~134--142.
\bibitem[LLS]{LLS} Li, C., Li, S.\ and Saito, K., Primitive forms via polyvector fields, arXiv:1311.1659.
\bibitem[Ma1]{M1} Malgrange, B., Polyn\^ome de Bernstein-Sato et cohomologie \'evanescente, Analysis and topology on singular spaces, II, III (Luminy, 1981), Ast\'erisque 101--102 (1983), 243--267.
\bibitem[Ma2]{M2} Malgrange, B., D\'eformations de syst\`emes diff\'erentiels et microdiff\'erentiels, Mathematics and physics, Progr.\ Math.\ 37, Birkh\"auser, Boston, MA, 1983, 353--379.
\bibitem[Ma3]{M3} Malgrange, B., Deformations of differential systems, II, J.\ Ramanujan Math.\ Soc.\ 1 (1986), 3--15.
\bibitem[Ph]{P} Pham, F., Singularit\'es des Syst\`emes Diff\'erentiels de Gauss-Manin, Prog.\ Math.\ 2, Birkh\"auser, Boston, MA, 1979.
\bibitem[Sab]{Sab} Sabbah,~C., Isomonodromic deformations and Frobenius manifolds. An introduction (Translated from the 2002 French edition), Universitext, Springer, London, 2007.
\bibitem[SaSa]{sasa} Sabbah,~C.\ and Saito, M., Kontsevich's conjecture on an algebraic formula for vanishing cycles of local systems, Algebraic Geometry 1 (2014), 107--130.
\bibitem[SK1]{SK1} Saito, K., Primitive forms for a universal unfolding of a function with an isolated critical point, J.\ Fac.\ Sci.\ Univ.\ Tokyo Sect.\ IA Math.\ 28 (1981), 775--792.
\bibitem[SK2]{SK2} Saito, K., Period mapping associated to a primitive form, Publ.\ Res.\ Inst.\ Math.\ Sci.\ 19 (1983), 1231--1264.
\bibitem[Sa1]{mhp} Saito, M., Modules de Hodge polarisables, Publ.\ RIMS, Kyoto Univ. 24 (1988), 849--995.
\bibitem[Sa2]{exp} Saito, M., Exponents and Newton polyhedra of isolated hypersurface singularities, Math.\ Ann.\ 281 (1988), 411--417.
\bibitem[Sa3]{bl} Saito, M., On the structure of Brieskorn lattice, Ann.\ Inst.\ Fourier 39 (1989), 27--72.
\bibitem[Sa4]{period} Saito, M., Period mapping via Brieskorn modules. Bull.\ Soc.\ Math.\ France 119 (1991), 141--171.
\bibitem[Sa5]{micro} Saito, M., On microlocal $b$-function, Bull.\ Soc.\ Math.\ France 122 (1994), 163--184.
\bibitem[SKK]{SKK} Sato, M., Kawai, T.\ and Kashiwara, M., Microfunctions and pseudo-differential equations, Lect. Notes in Math.\ 287, Springer, Berlin, 1973, pp.~265--529.
\bibitem[ScSt]{SS} Scherk, J.\ and Steenbrink, J.H.M.,
On the mixed Hodge structure on the cohomology of the Milnor fibre
Math.\ Ann.\ 271 (1985), 641--665.
\bibitem[St]{St} Steenbrink, J.H.M., Mixed Hodge structure on the vanishing cohomology, in Real and complex singularities, Sijthoff and Noordhoff, Alphen aan den Rijn, 1977, pp. 525--563.
\bibitem[Va]{V} Varchenko, A.N, Asymptotic Hodge structure in the vanishing cohomology, Math.\ USSR-Izv.\ 18 (1982), 469--512.
\end{thebibliography}
\end{document}